\documentclass{elsart}
\usepackage{amssymb, amsmath, amsfonts}
\usepackage[all]{xy}
\usepackage{bm}
\def\genfd{{\bm k}}
\def\id{{\rm id}}
\long\def\nodo#1{{}}
\def\gg{\mathfrak{g}}
\def\bbC{\mathbb{C}}
\def\CE{\mathbb{C}}
\def\PPartial#1{\frac{\partial}{\partial(\partial^{#1})}}
\def\Hom{\operatorname{Hom}}

\def\Der{\operatorname{Der}}
\def\Coder{\operatorname{Coder}}
\def\End{\operatorname{End}}
\def\hx{\hat{x}}
\def\hy{\hat{y}}
\def\ad{\operatorname{ad}}
\def\nxpoint{\refstepcounter{subsection}\makepoint{\thesubsection}}
\def\nxsubpoint{\refstepcounter{subsubsection}%
  \makepoint{\thesubsubsection}}
\def\refpoint#1{{\rm\textbf{\ref{#1}}}}
\def\makepoint#1{\medbreak\noindent{\bf #1. }}

\def\MR#1{} 
\def\hxi{\hat{x}_i}
\def\hxj{\hat{x}_j}
\def\hxk{\hat{x}_k}
\def\hf{\hat{f}}
\def\hpart{\hat\partial}
\begin{document}
\begin{frontmatter}
\title{Leibniz rules for enveloping algebras in symmetric ordering}
\author{Stjepan Meljanac, Zoran \v{S}koda}
\address{Theoretical Physics Division,
Institute Rudjer Bo\v{s}kovi\'{c}, Bijeni\v{c}ka cesta~54, P.O.Box 180,
HR-10002 Zagreb, Croatia
}
\ead{meljanac@irb.hr, zskoda@irb.hr}

\date{{}}

\maketitle
\begin{abstract}
Given a finite-dimensional Lie algebra, and a representation
by derivations on the completed symmetric algebra of its dual,
a number of interesting twisted constructions appear: certain
twisted Weyl algebras, deformed Leibniz rules, quantized ``star''
product. We first illuminate a number of interrelations between
these constructions and then proceed to study a special
case in certain precise sense corresponding
to the symmetric or Weyl ordering. This case has been known
earlier to be related to computations with Hausdorff series,
for example the expression for the star product is in such
terms. For the deformed Leibniz rule, hence a coproduct, we present
here a new nonsymmetric expression, which is then expanded
into a sum of expressions labelled by a class of planar trees,
and for a given tree evaluated by Feynman-like rules.
These expressions are filtered by a bidegree and we show recursion
formulas for the sums of expressions of a given bidegree, and
compare the recursions to recursions for Hausdorff series,
including the comparison of initial conditions. This way we
show a direct corespondence between the Hausdorff series and
the expression for twisted coproduct.
\end{abstract}
\end{frontmatter}

\section{The data defining the setup}

\nxpoint Fix a $n$-dimensional Lie algebra $\gg$ over a field $\genfd$.
The main message in our first several pages consists of 
the correspondences between several kinds of data:
\begin{itemize}
\item $\genfd$-linear maps $\bm\phi : \gg\to\Hom_\genfd(\gg,\hat{S}(\gg^*))$
\item $\genfd$-linear maps $\tilde{\bm\phi} : 
\gg\to\Der_\genfd(\hat{S}(\gg^*),\hat{S}(\gg^*))$
\item Matrices $(\phi^\alpha_\beta)_{\alpha,\beta = 1,\ldots, n}$ of elements 
$\phi^\alpha_\beta\in\hat{S}(\gg^*)$ satisfing the 
system of formal differential equations~(\ref{eq:phiderphi}).
\item Hopf actions of $U(\gg)$ on $\hat{S}(\gg^*)$.
\item Algebra homomorphisms $U(\gg)\to \hat{A}_{n,\genfd}$ 
(the codomain is the $n$-th Weyl algebra completed with respect to the powers
of $\partial^i$-s) which is of the form 
$\hx_\mu \mapsto \sum_{\alpha = 1}^n x_\alpha \phi^\alpha_\mu$ 
on a basis $\hx_1,\ldots,\hx_n$ of $\gg$, with 
$\phi^\alpha_\beta\in \hat{S}(\gg^*) $
\item Coalgebra isomorphisms $\xi : S(\gg)\to U(\gg)$ which are 
identity on $\gg \oplus\genfd = U^1(\gg)\subset U(\gg)$.
\end{itemize}
These correspondences are pretty easy to observe and 
the list can be meaningfully extended. On the other
hand, the many special cases of such data studied in (mainly
recent physics) literature, are treated with confusion about 
the definitions, nature and correspondences between these data and
related constructions. 

\nxsubpoint Note that the list can be meaningfully extended.
For example, there are popular ``ordering prescriptions'' which are 
various concrete ways determining the coalgebra isomorphism $\xi$ above
(as the isomorphism is trivial on generators in $\gg$ one needs to know
what to do with higher polynomials, hence ``ordering prescriptions'').
Another set of data, more loosely defined: 
if one extends $\xi$ to the completion (power series) then
$\xi$ can be evaluated on some interesting dense set, and exponentials
$\exp(ikx)$ are a good candidate and $\xi(\exp(ikx))$ is of the form
$\xi(\exp(iK(k)x))$ where $K : \genfd^n\to\genfd^n$ is a bijection
which is determined by $\phi$ and determines $\phi$ (however we do not
know {\it general} rule which bijections $K$ are admisible, though we
do have classification results for some very special $\gg$). 
For the correspondences to be bijections we need the assumption
that the maps $\bm\phi$ etc. are close to the ``unit'' case: for example 
in the case of $\phi^\alpha_\beta$ the unit case is $\delta^\alpha_\beta$
and the near by is in the sense of topology for the ring of
formal power series in $\partial$. Nonformal case of the correspondences
is interesting as well, but more difficult and we have 
no closed sufficiently general results of that form. 

\nxpoint {\bf (Morphism $\bm\phi$ and the equation it satisfies)}
Suppose we are given a Lie algebra $\gg$ and a
finite-dimensional vector space $V$ over a field $\genfd$.
By $\widehat{S(V)}$ or $\hat{S}(V)$ we will denote 
the completed symmetric algebra on $V$,
which may be viewed as a formal power series
in $m = \mathrm{dim}\,V$ variables.
Later we will set $V = \gg^*$,
but for the moment we consider the full generality.
Suppose we are also given a linear map
${\bm \phi} : \gg\to \Hom_\genfd(V,\widehat{S(V)})$.
We want to extend this
map to a $\genfd$-linear map into continuous derivations also denoted
$\tilde{\bm \phi} : \gg \to \Der_\genfd(\widehat{S(V)},\widehat{S(V)})$.
By the commutativity of $\widehat{S(V)}$, it must hold that
\begin{equation}\label{eq:onprod}
\tilde{\bm \phi}(\hat{x})(v_1\cdots v_n) = \sum_{i = 1}^n v_1\cdots v_{i-1}
 v_{i+1}\cdots v_n {\bm \phi}(\hat{x})(v_i).
\end{equation}
This formula is linear in all arguments and symmetric under their
permutations, hence by linearity in all arguments it defines a unique
extension of $\bm\phi(\hx)$ to a well-defined map
$\tilde{\bm\phi}(x)\in\Hom_\genfd(S(V),\widehat{S(V)})$.
It is straigthforward to check that $\tilde{\phi}(\hx)$ defined
via~(\ref{eq:onprod}) is indeed a derivation.
By abuse of notation, we will henceforth denote the extension
$\tilde{\bm\phi}$ also by $\bm\phi$.

Let $\partial^1,\ldots,\partial^m$ be a vector space basis of $V$.
Then, in terms of (algebraically defined) partial derivatives
$\PPartial{i}$, the condition~(\ref{eq:onprod}) generalizes to
the usual chain rule on $\widehat{S(V)}$
\begin{equation}
{\bm \phi}(\hat{x})(f)
= \sum_{i = 1}^m \PPartial{i}(f) {\bm\phi}(\hat{x})(\partial^i)
\end{equation}
Finally, we continuously extend $\bm\phi$ to a $\genfd$-linear map
$\bm \phi : \gg \to \operatorname{Der}_\genfd(\widehat{S(V)},\widehat{S(V)})$.

The enveloping algebra $U(\gg)$ is a Hopf algebra
with elements of $\gg\hookrightarrow U(\gg)$ being primitive.
If the linear map ${\bm \phi} : \gg\to \Der_\genfd(\widehat{S(V)})$ is
a homomorphism of Lie algebras, i.e.
\begin{equation}\label{eq:lieaction}
\bm\phi(\hx)\bm\phi(\hy)-\bm\phi(\hy)\bm\phi(\hx) - \bm\phi([\hx,\hy]) = 0,
\,\,\,\,\,\,\,\hx,\hy\in\gg,
\end{equation}
then $\bm\phi$ extends multiplicatively to a unique Hopf action of $U(\gg)$,
i.e. to a homomorphism
${\bm \phi} : U(\gg)\to \End_\genfd(\widehat{S(V)})$ satisfying
$\bm\phi(u)(fg) = m_{\hat{S}(V)}(\bm\phi\otimes\bm\phi)\Delta(u)(f\otimes g)=
\sum \bm\phi(u_{(1)})(f)\bm\phi(u_{(2)})(g)$,
for all $f,g \in \widehat{S(V)}$ and $u \in U(\gg)$,
where $m_{U(\gg)}$ is the multiplication map on $U(\gg)$.
From now on, let $\gg$ be finite-dimensional as well and
let $\hx_1,\ldots,\hx_n$ be a $\genfd$-basis of $\gg$,
and $\partial^1,\ldots,\partial^m$ the basis of $V$ (the Weyl algebras
usage later on for the case $V = \gg^*$ suggested the notation).
Denote
$$
\phi^\alpha_\beta = \phi^\alpha_\beta(\partial^1,\ldots,\partial^m)
:= \bm\phi(-\hx_\beta)(\partial^\alpha) \in \widehat{S(V)}.
$$
The formal power series
$\phi^\alpha_\beta =  \phi^\alpha_\beta(\partial^1,\ldots,\partial^m)$
has algebraically defined partial derivatives
$$
\PPartial{i} \phi^\alpha_\beta \in \widehat{S(V)}.
$$
Then $\bm\phi(\hxi)\bm\phi(\hxj)(\partial^k) =
\bm\phi(\hxi)(-\phi^k_j) = -\PPartial{l}(\phi^k_j)\bm\phi(\hxi)(\partial^l) =
- \PPartial{l}(\phi^k_j)\phi^l_i$.
Thus the condition~(\ref{eq:lieaction}) reads for $\hx = \hxi$ and $\hx =\hxj$
\begin{equation}\label{eq:phiderphi}
\phi^l_j \PPartial{l}(\phi^k_i) - \phi^l_i \PPartial{l}(\phi^k_j)
= C^s_{ij}\phi^k_s,
\end{equation}
where the structure constants $C^k_{ij}$ are defined by 
$[\hxi,\hxj] = C^k_{ij}\hxk$.
Consider the usual Weyl algebra $A_{n,\genfd}$ with generators
$x_1,\ldots,x_n,\partial^1,\ldots,\partial^n$, and its
completion $\hat{A}_{n,\genfd}$ with respect to the filtration 
by the degree of differential operator. Then the correspondence
$\hx_i \mapsto \hx_i^\phi := \sum_{j=1}^m x_j \phi^j_i$
extends to an algebra homomorphism
$()^\phi : U(\gg)\to \hat{A}_{n,\genfd}$
iff~(\ref{eq:phiderphi}) holds.

\nxsubpoint 
Many solutions for $\bm\phi$ satisfying~(\ref{eq:phiderphi}) and hence
homomorphisms $\bm\phi$, injective or not, for particular $\genfd$ and
particular $\genfd$-Lie algebras, with $m$ equal $n$ or not, exist.
For example, in the article of Berceanu~\cite{Berceanu},
such realizations with $\phi$ faithful have been found for $\gg$ semisimple
over $\genfd = \CE$. In that case, $\mathrm{dim}\, V < \mathrm{dim}\,\gg$
and $\mathrm{dim}\,V$ may be calculated in
terms of the combinatorics of root systems.

\nxpoint A universal formula for a "symmetric" solution
to~(\ref{eq:phiderphi}) has been found~(\cite{ldWeyl}),
for any {\it ring} $\genfd\supset \mathbb Q$,
$\gg$ finite rank free module over $\genfd$
and $V = \gg^*$ (in particular, $m = n$), where
$\bm\phi$ is a monomorphism. See also Section~\ref{secsymm}.

\nxpoint {\bf (Hopf algebras)} All bialgebras in the article will be
associative, coassociative, with unit map $\eta$ and counit
$\epsilon$, without gradings.
Hopf algebras will be bialgebras with an antipode and the standard
Sweedler notation for th coproduct
$\Delta(h) = \sum h_{(1)}\otimes h_{(2)}$ is few times used,
with or without the summation sign. Recall that the elements $h\in H$
such that $\Delta(h) = 1\otimes h + h\otimes 1$ are called primitive.

\nxpoint {\bf (Smash product algebras)}
Given any Hopf algebra $H$ and a, say left, Hopf action of $H$ on
algebra $\mathcal S$, $h\otimes s \mapsto h\triangleright s$,
one forms a crossed product algebra (in Hopf literature
"smash product") $\mathcal S\sharp H$. As a vector space, it is simply
the tensor product vector space $S\otimes \mathcal H$ and
the associative product is given by
$$
(s\otimes h)(s'\otimes h') =
\sum s (h_{(1)}\triangleright s')\otimes h_{(2)}h'.
$$
The canonical embeddings
$S\hookrightarrow S\sharp H$ and $H\hookrightarrow S\sharp H$
will be considered identifications, and one usually omits the tensor sign
because $s\otimes h = sh$ with respect to these embeddings
and the product in $\mathcal S\sharp H$.
Then $h\triangleright s = \sum h_{(1)} s S_H(h_{(2)})$
where $S_H : H\to H$ is the antipode.
Furthermore, the rule
\begin{equation}\label{eq:smashacts}
(s\sharp h)\triangleright s' := s (h\triangleright s')
\end{equation}
defines an action of $\mathcal S \sharp H$ on $\mathcal S$.

Analogously, for any {\it right} action of $H$ on $\mathcal S$
one defines the crossed product denoted by $H\sharp\mathcal S$,
whose underlying vector space is $H\otimes\mathcal S$.
If the antipode $S_H : H\to H^{\mathrm{op}}$ is bijective,
there is a bijective correspondence between the left and right actions
(namely, composing with $S_H$) and the crossed products
for the two corresponding (left and right) actions are canonically
isomorphic and we often identify them throughout the article.

\nxpoint {\bf ($(\gg, {\bm \phi})$-deformed Weyl algebras.)}
Regarding that for any $\gg$, $V$ and $\bm\phi$
such that~(\ref{eq:phiderphi}) holds,
the action of $U(\gg)$ on $\widehat{S(V)}$ is a Hopf action,
we may define the smash product algebra
$$
A_{\gg,\phi} := \widehat{S(V)}\sharp\, U(\gg) =
\widehat{S(V)}\sharp_\phi U(\gg),
$$
where the action $u\triangleright v := \bm\phi(u)(v)$
is uniquely determined by the values $\bm\phi(-\hx_i)(\partial^j) = \phi^i_j$
as explained above.
The rule~(\ref{eq:smashacts}) specializes to
a (dual) "natural" action of $A_{\gg,\phi}$ on $\widehat{S(V)}$.
In particular, if $\gg = \mathfrak a$ is an abelian Lie algebra,
$V = \gg^*$ and $\bm\phi$ is given simply by
the bilinear pairing $\bm\phi(\hxi)(\partial^j) = \delta^j_i$, then
$A_{\gg,\phi}$ is isomorphic to the usual
(semi)completed Weyl algebra $\widehat{A}_{n,\genfd}$ and
the action is the usual action of $S({\mathfrak a})$ on
$\widehat{S({\mathfrak a}^*)}$.

\nxpoint \label{s:standardphi} From now on we suppose

(i) $\bm\phi :\gg\to\Der(\widehat{S(\gg^*)})$ is a homomorphism
of Lie algebras

(ii) the matrix $\phi$ (not bold) with entries
$\phi^i_j := \phi(-\hxi)(\partial^j)$ has the unit matrix as
its constant term, i.e. $\phi^i_j = \delta^i_j + O(\partial)$.

\nxpoint \label{s:identification}
Under the assumptions from~\refpoint{s:standardphi},
$\phi$ is invertible as a matrix over the formal power series
ring $\genfd[[\partial^1,\ldots,\partial^n]]$ and the homomorphism
$U(\gg)\sharp\widehat{S(\gg^*)}\cong S(\gg)\sharp\widehat{S(\gg^*)}$
given on generators by
$$
\hat{x}_\alpha\mapsto x_\beta \phi^\beta_\alpha,\,\,\,\,\,
\partial^\mu\mapsto\partial^\mu
$$
is an isomorphism. Hence the (one-sidedly) completed deformed
and undeformed Weyl algebras are isomorphic via a nontrivial map
and we often identify them when doing calculations.

\nxpoint \label{s:actionUonS}
This isomorphism enables us to consider the homomorphism
$$
()^\phi : U(\gg)\hookrightarrow U(\gg)\sharp\widehat{S(\gg^*)}\cong
S(\gg)\sharp \widehat{S(\gg^*)} = \hat{A}_{n,\genfd}
$$
which agrees with the unique homomorphism $U(\gg)\to \hat{A}_{n,\genfd}$
extending the rule
$$
\hat{x}_\alpha \mapsto \hat{x}_\alpha^\phi :=
x_\beta\phi^\beta_\alpha
\in \hat{A}_{n,\genfd}
$$
Furthermore, we may identify
$S(\gg)\sharp \widehat{S(\gg^*)}\cong \Hom_\genfd(S(\gg),S(\gg))$.
Here $\phi^\beta_\alpha = \phi^\beta_\alpha(\partial^1,\ldots,\partial^n)$
is understood as an element of the completed Weyl algebra
$\widehat{A}_{n,\genfd} \cong S(\gg)\sharp \widehat{S(\gg^*)}$
acting in the usual way (as differential operator,
this one with constant coefficients) on the polynomial algebra.
Therefore we obtained an action, depending on $\phi$, of $U(\gg)$
on $S(\gg)$.

\nxpoint {\bf Lemma.} {\it
Let $\chi\in \genfd [[\partial^1,\ldots,\partial^n]]$, then
$x_\sigma \chi$ is a coderivation of the polynomial algebra
$P = \genfd[ x_1,\ldots,x_n ]$. In other words,
}
\begin{equation}
(x_\sigma \chi\otimes \id + \id\otimes x_\sigma \chi)
(\Delta_{P}(f)) = \Delta_{P}(x_\sigma \chi(f)),
\,\,\,\,\,\forall f\in P.
\end{equation}
{\it Proof.} By linearity it is enough to prove it for
$f$ which are monomial.
We prove this by induction on
the sum of the polynomial degree of $f$
and the order of differential operator $\Xi$.
The identity is clearly true if either the degree of $f$
or order of $\xi$ is $0$. Regarding that $f$ is monomial
it is of the form $x_\gamma g$ where $g$ is some monomial of
a lower order. We identify $\id$ with $1$ in Weyl algebra
and $x_\mu$ with multiplication with $x_\mu$.
For step of induction we want to prove that
$$
(x_\sigma\xi\otimes 1 + 1\otimes x_\sigma\xi)\Delta(x_\gamma g)
= \Delta(x_\sigma\xi(x_\gamma g))
$$
provided this is true for $\xi$ of lower order or $x_\gamma g$ replaced
by $g$ what is of lower degree. Using the fact that $\Delta$ is a
 homomorphism of algebras and that $x_\gamma$ is primitive,
we rewrite this equality using commutators:
$$\begin{array}{lcl}
(x_\sigma [\xi,x_\gamma]\otimes 1 &+& 1\otimes x_\sigma [\xi,x_\gamma])
\Delta(g) + \\
& + & (x_\gamma\otimes 1 + 1\otimes x_\gamma)
(x_\sigma \xi\otimes 1 + 1\otimes x_\sigma\xi)\Delta(g)\\
 && =  \Delta(x_\sigma [\xi,x_\gamma] (g)) +
(x_\gamma\otimes 1 + 1\otimes x_\gamma)
\Delta(x_\sigma\xi(g))
\end{array}$$
and recall that $[\xi,x_\gamma]$ is of lower order. This equality
holds because it is a sum of two equations which hold by the assumption
of the induction. Q.E.D.

{\bf Corollary.} {\it
The action from \refpoint{s:actionUonS} restricted on $\gg$
is an action by coderivations with respect to the standard coalgebra
structure on $S(\gg)$:
}
\begin{equation}\label{eq:actcoderiv}
(x_\beta \phi^\beta_\alpha\otimes \id + \id\otimes x_\beta \phi^\beta_\alpha)
(\Delta_{S(\gg)}(f)) = \Delta_{S(\gg)}(x_\beta \phi^\beta_\alpha(f)),
\,\,\,\,\,\forall f\in S(\gg).
\end{equation}

\nxpoint \label{s:xi}
For us it is important to consider the special
case of the action of $U(\gg)$ on $S(\gg)$
from~\refpoint{s:actionUonS}, when $f$ is $1$ (action on ``vacuum'').
In the following, recall the notation 
$()^\phi : u \mapsto u^\phi\in\hat{A}_{n\genfd}$ 
from \refpoint{s:actionUonS} and the evaluation on $1\in S(\gg)$ 
in $\hat{u}^\phi(1)$ is in the sense of natural action of $\hat{A}_{n\genfd}$ 
on $\hat{S}(\gg)$. 

{\bf Proposition.} {\it The rule $\xi^{-1} :\hat{u}\mapsto \hat{u}^\phi(1)$
for $u \in U(\gg)$ is an isomorphism of coalgebras, which restricts to
the identity on $\genfd\oplus\gg$. 
}

Of course, the inverse of $\xi^{-1}$ will be some isomorphism of coalgebras
$\xi : S(\gg)\to U(\gg)$. Conversely,  every isomorphism
of coalgebras $\xi : S(\gg)\to U(\gg)$ which is identity on $\genfd\oplus\gg$,
defines a map $D^T : \gg\to\Coder(S(\gg))$ into coderivations by
$D^T_x (f) =  D^T(x)(f) = \xi^{-1}\left(\xi(x)\cdot_{U(\gg)}\xi(f)\right)$.
The dual map $D_x : \widehat{S(\gg^*)}\to\widehat{S(\gg^*)}$
is a continuous derivation, and one has
$D^T_x(f) = -\sum_\alpha x_\alpha D_x(\partial^\alpha)(f)$ where the action on
the left is the usual action as differential operator.
Here $\sum_\alpha x_\alpha \otimes \partial^\alpha\in\gg\otimes\gg^*$ is
the ``canonical element'' (the image of $\id_\gg$ under the isomorphism
$\Hom_\genfd(\gg,\gg)\to\gg\otimes\gg^*$).
Thus one defines a Lie homomorphism
$\bm\phi : \gg\to\Der(\widehat{S(\gg^*)},\widehat{S(\gg^*)})$ by $x\mapsto D_x$
such that $\phi^i_j = D_{x_j}(\partial^i)$ and
$\phi^i_j = \delta^i_j + O(\partial)$.

\nxpoint {\bf (Star product)} We saw in \refpoint{s:xi}
that giving the Lie homomorphism $\bm\phi$ for which the matrix
$\phi(-\hxi)(\partial^j) = \delta^i_j + O(\partial)$ is equivalent to
giving a coalgebra isomorphism $\xi : S(\gg)\to U(\gg)$
which is identity when restricted to $\genfd\oplus\gg$.
This isomorphism helps us define the star product
\begin{equation}\label{eq:star}
\star : S(\gg)\otimes S(\gg)\to S(\gg),\,\,\,\,\,
f\star g := \xi^{-1}(\xi(f)\cdot_{U(\gg)}\xi(g)).
\end{equation}

\nxsubpoint One should note that in literature related to the 
representation theory~(\cite{Olshanski}) and
the deformation quantization~(\cite{KontsDefPoiss})
usually some other isomorphisms $S(\gg)\to U(\gg)$
are important, which do not respect the coalgebra structure,
but do have some other favorable properties. 
Our constructions below, however, essentially use the
compatibility with coalgebra structure. 

\nxpoint \label{s:actionSonU}
Given a Hopf algebra $H$ acting by a right Hopf action
on an algebra $\mathcal S$ and a
homomorphism of unital algebras $\epsilon^S :\mathcal S\to\genfd$, one
defines a $\genfd$-linear map
$(H\sharp\mathcal S)\otimes H\to H$
as the following composition:
$$
(H\sharp\mathcal S)\otimes H\hookrightarrow
(H\sharp\mathcal S)\otimes (H\sharp\mathcal S)
\stackrel{m_{H\sharp\mathcal S}}\longrightarrow H\sharp\mathcal S
\stackrel{H\sharp\epsilon^S}\longrightarrow H\otimes\genfd\cong H.
$$
This map is a left{\it action} of
the smash product algebra $H\sharp\mathcal S$ on $H$.
Algebra embedding $S\hookrightarrow H\sharp\mathcal S$, $s\mapsto 1\otimes s$,
gives rise to the restriction of the above action to
a left action $\mathcal S\otimes H\to H$. If the antipode
$S_H: H\to H^{\mathrm{op}}$ is an isomorphism, the corresponding
representation $\rho :\mathcal S\to\End_\genfd(H)$ is faithful.

\section{Deformed derivatives}

\nxpoint We now define the deformed derivatives in several ways.
Physicist who view $U(\gg)$ as some sort of algebra of
functions on ``Lie type noncommutative space'' consider the
deformed derivatives (which mutually commute) as a sensible 
choice of a basis of the tangent space to this noncommutative space
~(\cite{DimGauge,HallidaySzabo,covKappadef,kappaind,MS}). 

\nxpoint \label{s:hpart}
In the case when $\mathcal S = \widehat{S(\gg^*)}$, $H = U(\gg)$
and the Hopf action is induced by
$\bm\phi : U(\gg)\to\Der(\widehat{S(\gg^*)},\widehat{S(\gg^*)})$,
the action $\rho_\phi : \mathcal S\to\End_\genfd(H)$
from~\refpoint{s:actionSonU}
may be alternatively described in terms of 
the values on the standard algebra generators
$\hpart^\mu = \rho_\phi(\partial^\mu)\in\End_\genfd(U(\gg))$,
$\mu = 1,\ldots,n$. 

We describe the action of $\hpart^\mu$ on $U(\gg)$
inductively on the order of monomials in $U(\gg)$. 
First of all, $\hat\partial^\mu(1) = 0$ and
$\hat\partial^\mu(\hat{x}_\nu)=\delta^\mu_\nu$.
Then suppose $\hpart^\mu$ is defined on monomials of order up to $n$.
Then any monomial of order $n+1$ is of the form $\hx_\nu \hf$ where
$\hpart(\hf)$ is already defined. We set
$$
\hpart^\mu(\hx_\nu \hf) := [\hpart^\mu, \hx_\nu](\hf) +
\hx_\nu \hpart^\mu(\hf)
:= \phi^\mu_\nu(\hf) + \hx_\nu \hpart^\mu(\hf),
$$
where $\phi^\mu_\nu = \phi^\mu_\nu(\hpart)$
(we can substitute $\hpart$
because $S(\gg^*)$ is a free commutative algebra
and $\hat\partial^\mu$ mutually
commute as it may be shown a posteriori).
$\hpart$ is well defined on $S(\gg^*)$ (hence by continuity on
$\widehat{S(\gg^*)}$), namely it is obviously well defined linear operator
from the free algebra on abstract variables $\hx_\alpha$ to $U(\gg)$,
and if one takes the generators of the defining ideal of the enveloping algebra
$i_{\nu_1\nu_2} = \hat{x}_{\nu_1} \hat{x}_{\nu_2} -
\hat{x}_{\nu_2}\hat{x}_{\nu_1} - C^\alpha_{\nu_1\nu_2} \hat{x}_\alpha$
then, applying our inductive rules for every of the three monomials on RHS,
we conclude that for every $\hf\in U(\gg)$,
$$\begin{array}{lcl}
\hat\partial^\gamma(i_{\nu_1\nu_2}\hf) & = &
\phi^\gamma_{\nu_1}(\hx_{\nu_2}\hf) +
\hx_{\nu_1} \hat\partial^\gamma\hx_{\nu_2}(\hf)
- \phi^\gamma_{\nu_2}(\hx_{\nu_1}\hf) - \\
&& \,\,\,\,\,\,\, - \hx_{\nu_2} \hat\partial^\gamma\hx_{\nu_1}(\hf)
- C^\alpha_{\nu_1\nu_2} \phi^\gamma_\alpha (\hf)
- C^\alpha_{\nu_1\nu_2} \hx_\alpha \hat\partial^\gamma(f) \\
& = &
\PPartial{\nu_2}(\phi^\gamma_{\nu_1})(\hf) +
\hx_{\nu_2} \phi^\gamma_{\nu_1}(\hf) + \hx_{\nu_1}\phi^\gamma_{\nu_2}(\hf)
+ \hx_{\nu_1} \hx_{\nu_2} \hpart^\gamma (\hf)\\
&&\,\,\,\,\,\,\,
-  \left(\PPartial{\nu_1}(\phi^\gamma_{\nu_2})(\hf) +
\hx_{\nu_1} \phi^\gamma_{\nu_2}(\hf) + \hx_{\nu_2}\phi^\gamma_{\nu_1}(\hf)
+ \hx_{\nu_1}\hx_{\nu_2}\hpart^\mu(\hf)\right)\\
&& \,\,\,\,\,\,\, - C^\alpha_{\nu_1\nu_2} \phi^\gamma_\alpha (\hf)
- C^\alpha_{\nu_1\nu_2} \hx_\alpha \hat\partial^\gamma(\hf)\\
&=& (\PPartial{\nu_2}(\phi^\gamma_{\nu_1}) -
\PPartial{\nu_1}(\phi^\gamma_{\nu_2}) -
C^\alpha_{\nu_1\nu_2} \phi^\gamma_\alpha )(\hf)
\end{array}$$
The injectivity of $\rho$ implies that
$\hat\partial^\gamma(i_{\mu\nu}\hf) = 0$  for every $\hf$ iff
the operator in the brackets on RHS vanishes,
what amount to our main assumption~(\ref{eq:phiderphi}).
It is trivial that $\hat\partial(\hf i_{\mu\nu}) = 0$
as well, namely this is sufficient to check for monomial $\hf$, but this is
$\hf \hat\partial (i_{\mu\nu}) + [\hat\partial,\hf](i_{\mu\nu})$.
We already know that the first summand is zero.
The commutator in the second summand
is some polynomial in $\hpart$-s, hence it is clearly zero modulo
$i_{\mu\nu}$ by induction on the degree of monomials and linearity.

Notice for the classical case of the abelian Lie algebra,
that $[\partial, \hat{f}] = \hat{\partial}(\hat{f})$,
while this is not true in general (the equality always makes sense:
LHS is the bracket $\partial \hat{f}-\hat{f}\partial$
in the smash product $\widehat{S(\gg^*)}\sharp U(\gg)$, while RHS is in
$U(\gg)\hookrightarrow \widehat{S(\gg^*)}\sharp U(\gg)$).

\nxpoint One can alternatively describe operators
$\hpart^\mu = \rho_\phi(\partial^\mu)$
from \refpoint{s:hpart} by the formula
$$
\hpart^\mu(\xi(f)) = \xi (\partial^\mu(f)), \,\,\,\,\,\,f\in S(\gg),
$$
where $\xi = \xi_\phi$ is described in~\refpoint{s:xi}.
Therefore also $\xi^{-1}\hpart^\mu = \partial^\mu\xi^{-1}$.
It is straightfoward to check that this description agrees with
the inductive description of $\hpart^\mu$ in~\refpoint{s:xi}.
It is of course convenient to have such an invariant description.
Moreover, for any $P \in \hat{A}_{n\genfd}$, 
understood in a usual way as an operator $S(\gg)\to S(\gg)$, 
we form $\hat{P} : U(\gg)\to U(\gg)$ 
by the same transport rule $\hat{P}(\xi(f)) = \xi(P(f))$.

The deformed coproduct $\Delta(\hpart^\mu) = \sum \hpart^\mu_{(1)}
\otimes\hpart^\mu_{(2)}$ is defined by
$\hpart(u\cdot_{U(\gg)} v) =
\sum \hpart^\mu_{(1)}(u)\cdot_{U(\gg)} \hpart^\mu_{(2)}(v)$
for $u,v\in U(\gg)$. This is
equivalent to the ``deformed Leibniz rule'', 
popular in some physics works:
$$ \partial^\mu (f\star g) =
\sum_i \partial^\mu_{(1)}f\star \partial^\mu_{(2)}g,\,\,\,\,\,\,\,
\,\,\,\,f,g \in S(\gg),$$
as the following calculation shows:
$\partial^\mu(f\star g) =
\partial^\mu(\xi^{-1}(\xi(f)\cdot_{U(\gg)}\xi(g)) =
\xi^{-1}(\hpart^\mu(\xi(f)\cdot_{U(\gg)}\xi(g))) =
\xi^{-1}(\hpart^\mu_{(1)}\xi(f)\cdot_{U(\gg)}\hpart^\mu_{(2)}\xi(g))
= \xi^{-1}(\xi\hpart^\mu_{(1)}(f)\cdot_{U(\gg)}\xi\hpart^\mu_{(2)}(g))
= \partial^\mu_{(1)}(f)\star\partial^\mu_{(2)}(g)$.

From now on, when commuting with elements in
$U(\gg)\hookrightarrow A_{\gg,\phi}$ we will by
$\partial\in S(\gg^*)$ mean $\hpart\in S(\gg^*)\hookrightarrow A_{\gg,\phi}$
and so on -- this amounts to
the identification of the deformed and undeformed Weyl algebras,
cf.~\refpoint{s:identification}. The main purpose of this article
is describing more concretely this deformed coproduct.

\nxpoint This coproduct is related to
but different from the dual coproduct
$(S(\gg))^*\cong\widehat{S(\gg^*)}\stackrel{\Delta'}\to\widehat{S(\gg^*)}
\widehat\otimes\widehat{S(\gg^*)}$ to the star product~(\ref{eq:star}).
The defining property of $\Delta'$ is
$\langle u_{1'},f\rangle\langle u_{2'},g \rangle\equiv\langle\Delta'(u),
f\otimes g\rangle = \langle u, f\star g\rangle$ for
$u\in S(\gg)^*\cong \widehat{S(\gg^*)}$, $f, g\in S(\gg)$.

The coorespondence $P\mapsto (f\mapsto P(f)(0))$ is the linear
isomorphism from the space of derivations of $S(\gg)$ to the space of linear
functionals $S(\gg)^*$.  Evaluating at zero the $n$-th partial
derivative is the same as evaluating the product of first partial derivatives
except that one has to adjust the factor of $n!$ what amounts to
a different pairing between the graded components $S^n(\gg)$ and $S^n(\gg^*)$
(i.e. a different identification $S^n(\gg^*)\cong S^n(\gg)^*$).

\nxpoint {\bf Lemma.} {\it If $\hat{a} = a^\alpha \hat{x}_\alpha$
and $\hat{f}\in U(\gg)$ then}
\begin{equation}\label{eq:dxn}
\hpart^\mu (\hat{a}^p \hat{f}) =
\sum_{k = 0}^{p-1} {n\choose k} a^{\alpha_1} a^{\alpha_2}\cdots a^{\alpha_k}
\hat{a}^{p-k} [[[\hpart^\mu,\hx_{\alpha_1}],\ldots,\hx_{\alpha_k}] (\hat{f})
\end{equation}
{\it Proof.} This is a tautology for $p = 0$. Suppose it holds for all $p$
up to some $p_0$, and for all $\hat{f}$.
Then set $\hat{g} = \hat{a}\hat{f} = a^\alpha \hat{x}_\alpha$.
Then $\hpart^\mu(\hat{a}^{p_0 +1}\hat{f}) = $
and we can apply~(\ref{eq:dxn}) to $\hpart^\mu(\hat{a}^{p_0}\hat{g})$. Now
$$\begin{array}{lcl}
[[[\hpart^\mu,\hx_{\alpha_1}],\ldots],\hx_{\alpha_k}] (\hat{g})
&=&  a^{\alpha_k}[[[\hpart^\mu,\hx_{\alpha_1}],\ldots],\hx_{\alpha_k}]
(\hx_{\alpha_k}\hat{g})\\
&=&\hat{a}[[[\hpart^\mu,\hx_{\alpha_1}],\ldots],\hx_{\alpha_k}] (\hat{f})
\,+ \\&&\,\,\,\,\,\,\,+ \,\,a^{\alpha_{k+1}}
[[[[\hpart^\mu,\hx_{\alpha_1}],\ldots],\hx_{\alpha_k}],\hx_{\alpha_{k+1}}]
(\hat{f}).
\end{array}$$
Collecting the terms and the Pascal triangle identity complete the
induction step.

\section{Symmetric ordering}\label{secsymm}

\nxpoint \label{s:symmord}
Given a basis $\hx_1,\ldots,\hx_n$
in a Lie algebra $\gg$, and structure constants
defined by $[\hx_i,\hx_j] = C^k_{ij} \hx_k$,
denote by $\bbC$ the matrix with entries in $A_{n,\genfd}$
whose $(i,j)$-th entry is
$$
\bbC^i_j = C^i_{jk} \partial^k
$$
In~\cite{ldWeyl} we have shown that if $\xi:S(\gg)\to U(\gg)$
is the coexponential map then the corresponding $\bm\phi$ is determined by
$$
\bm\phi(-\hx_\beta)(\partial^\alpha) =
\phi^\alpha_{\beta} = \sum_{N=s}^\infty
(-1)^N \frac{B_N}{N!} (\CE^N)^\alpha_{\beta}
$$
where $B_N$ are the Bernoulli numbers. {\it For the reason of
structure of the coexponential map, from now on we will
say that this is the case of} {\bf symmetric ordering}.
It has the property that $\xi^{-1}(\exp(a^\alpha\hx_\alpha)) =
\exp(a^\alpha x_\alpha)$. In fact there is a bit more general fact,
which we will show in~\cite{snote}:

\nxpoint {\bf (Symmetric case; only tensorial form used)}
Given $C^k_{ij}, \bbC$ as above let $U$ be {\it any} subalgebra
of $A_{n,\genfd}[[t]]$ (a priori not necessarily isomorphic to $U(\gg)$)
generated by $n$ generators $X_1,\ldots,X_n$
which satisfies the following two conditions

(i) the mapping $x_{\alpha_1}\cdots x_{\alpha_k}\mapsto
\frac{1}{|\Sigma(k)| !}\sum_{\sigma \in \Sigma(k)}
X_{\alpha_{\sigma 1}}\cdots X_{\alpha_{\sigma k}}$  extends to
a onto map $\xi : \genfd[x_1,\ldots,x_n]\to U$

(ii) $X_i = \sum_{N=0}^\infty A_N x_\alpha (\bbC^N)^\alpha_i$,
where $A_N \in \genfd$ for all $N>0$
are arbitrary, $A_0 = 1$ and where the
summation over $\alpha$ is understood.
We will denote $\phi = \sum_{N=0}^\infty A_N \bbC^N$, hence
$X_i = x_\alpha \phi^\alpha_i$.

Then the following theorem holds

\nxsubpoint {\bf Theorem.} {\it
Let $\theta : U \to \genfd[x_1,\ldots,x_n]$
be defined as
$$
\theta(P) = P(1).
$$
where $P(1)$ is evaluated in the sense of
the natural action of $A_{n,\genfd}[[t]]$ on $\genfd[x_1,\ldots,x_n][[t]]$.
Then $\theta\circ\xi = \id$. In particular, $\xi$ is then
injective, hence by (i) an isomorphism of vector spaces.
}

\nxpoint For general $\phi$,  
$[\hat\partial^\mu,\hx_\alpha] = \phi^\mu_\alpha$,
$[\hat\partial^\mu,\hx_\alpha](1) = \delta^\mu_\alpha$, and

$$[[\hat\partial^\mu,\hx_\alpha],\hx_\beta] =
\frac{\partial}{\partial(\partial^\rho)}
(\phi^\mu_\alpha)\phi^\rho_\beta = \phi^\mu_{\alpha,\rho}\phi^\rho_\beta
$$
holds. 
In the case of the symmetric ordering (cf.~\refpoint{s:symmord}), 
that is when $\xi$ is the coexponential map~(\cite{ldWeyl}), also
$$
\phi^\mu_{\alpha,\rho}\phi^\rho_\beta (1) =
\frac{1}{2} C^\mu_{\alpha\beta}
$$
and the higher order terms are not so easy to evaluate at 1 in 
general in a closed form (this involves identities 
between different tensors in $C$-s, what is combinatorially involved,
hence one should probably  handled it using tree calculus).

\nxpoint (Notation: subscripts after comma for derivatives.)
Given $\phi^\alpha_\beta\in \widehat{S(\gg^*)}$ as above, denote
$$
\phi^\alpha_{\beta,\rho_1\rho_2\ldots\rho_k} :=
\frac{\partial}{\partial(\partial_{\rho_k})}\ldots
\frac{\partial}{\partial(\partial_{\rho_2})}
\frac{\partial}{\partial(\partial_{\rho_1})} \phi^\alpha_\beta
$$
and we use the extension of this notation
to more complicated expressions, e.g.
$(ab)_{,\rho} = a_{,\rho} b + a b_{,\rho}$
is the derivative of the product $ab$ with respect to $\partial_\rho$.

\nxpoint {\bf Lemma.} {\it
Let $\hx_1,\ldots,\hx_n$ be a basis of $\gg$.
For any $\phi$ as above,
}
\begin{equation}\label{eq:phiphibr}
[\ldots[[\hpart^\mu,\hx_{\alpha_1}],\hx_{\alpha_2}],\ldots,
\hx_{\alpha_k}]
= (\ldots((\phi^\mu_{\alpha_1,\rho_1}\phi^{\rho_1}_{\alpha_2})_{,\rho_2}
\phi^{\rho_2}_{\alpha_3})_{,\rho_3}\ldots )_{,\rho_{k-1}}
\phi^{\rho_{k-1}}_{\alpha_k}
\end{equation}

The proof is an obvious induction, using the chain rule.

\nxsubpoint 
Using the Leibniz rule we can rewrite the formula~(\ref{eq:phiphibr})
as a sum of terms for which every derivative operator
$\frac{\partial}{\partial(\partial_\rho)}$
is applied only to a single $\phi$-series, rather than to products.
Indeed, it is clear that $\frac{\partial}{\partial(\partial_{\rho_1})}$
applies only to $\phi^\mu_{\alpha_1}$, then
$\frac{\partial}{\partial(\partial_{\rho_2})}$
applies either to $\phi^\mu_{\alpha_1}$ or
$\phi^{\rho_1}_{\alpha_2}$, and in general,
$\frac{\partial}{\partial(\partial_{\rho_s})}$
applies to $\phi^{\rho_{p-1}}_{\alpha_p}$ where $1\leq p \leq s$
and $\rho_0 := \mu$. This means that we have $(k-1)!$ summands.
For example for $k = 4$ we have 6 summands:
$$\begin{array}{l}
\phi^\mu_{\alpha_1,\rho_1} \phi^{\rho_1}_{\alpha_2,\rho_2}
\phi^{\rho_2}_{\alpha_3,\rho_3}\phi^{\rho_3}_{\alpha_4} +
\phi^\mu_{\alpha_1,\rho_1} \phi^{\rho_1}_{\alpha_2,\rho_2\rho_3}
\phi^{\rho_2}_{\alpha_3}\phi^{\rho_3}_{\alpha_4} +
\phi^\mu_{\alpha_1,\rho_1\rho_3} \phi^{\rho_1}_{\alpha_2,\rho_2}
\phi^{\rho_2}_{\alpha_3}\phi^{\rho_3}_{\alpha_4}
\\ +
\phi^\mu_{\alpha_1,\rho_1\rho_2} \phi^{\rho_1}_{\alpha_2}
\phi^{\rho_2}_{\alpha_3,\rho_3}\phi^{\rho_3}_{\alpha_4} +
\phi^\mu_{\alpha_1,\rho_1\rho_2} \phi^{\rho_1}_{\alpha_2,\rho_3}
\phi^{\rho_2}_{\alpha_3}\phi^{\rho_3}_{\alpha_4} +
\phi^\mu_{\alpha_1,\rho_1\rho_2\rho_3}
\phi^{\rho_1}_{\alpha_2} \phi^{\rho_2}_{\alpha_3}\phi^{\rho_3}_{\alpha_4}
\end{array}$$
I will call this expansion ``expansion 1''.

\nxpoint We now specialize to the case of the series corresponding to
the symmetric ordering
$$
\phi^\alpha_{\beta,\rho_1,\ldots,\rho_s} = \sum_{N=s}^\infty
(-1)^N \frac{B_N}{N!} (\CE^N)^\alpha_{\beta,\rho_1\ldots\rho_s}
$$
The sum over $N\geq k$ for each $\phi$ in the form of expansion 1,
will be called expansion 2.
By applying the Leibniz rule again, we notice that
$
(\CE^N)_{,\rho_1\ldots\rho_s}
$
is a sum of $N!/(N-s)!$ summands, each of which is monomial
which is a product of $N-s$ $\CE$-s and $s$ $C$-s. This is
the expansion 3.
Performing consequently expansions 1,2 and 3,
the commutator in~(\ref{eq:phiphibr}) becomes a
multiple sum of terms which
are labelled by certain class
of attributed planar trees and each summand is certain contraction
of several $\CE$-tensors and several $C$-tensors with $k+1$ external
indices $\mu,\alpha_1,\ldots,\alpha_k$, and with some pre-factor
involving (products of) Bernoulli numbers and factorials.
To describe the details, we introduce several ``classes'' of
planar rooted trees and their ``semantics''.

\section{Tree calculus for symmetric ordering}

\nxpoint {\it Class $\mathcal T$ consists of all planar rooted trees
with two kinds of nodes, white and black, where black nodes may only be leaves.
}

We will draw the trees in $\mathcal T$ with the root on the top.
'Planar' implies that the (left to right) order of child branches
of every node matters. If $t\in\mathcal T$,
then $w(t)\geq 0$ and $b(t)\geq 0$ are
the number of white and black nodes in $t$ respectively. Class $\mathcal T$
is graded in obvious way $\mathcal T = \coprod_{P=1}^\infty \mathcal T_P$
by the total number of nodes $P$, and bigraded by the numbers $b$ and $w$ of
black and white nodes: $\mathcal T = \coprod_{w+b>0} \mathcal T_{w,b}$.
Clearly $\mathcal T_P = \coprod_{w+b = P} T_{w,b}$.

Class $\mathcal T^{\mathrm ord}$ consists of pairs $(t,l)$ where
$t\in\mathcal T$ and $l$ is a numeration (with values $1,\ldots,w$)
on the set of white nodes of $t$ which is descending in the sense that
white children nodes are always assigned
greater values than their parent nodes.
Let $\mathcal T^{\mathrm ord}_P$ and $\mathcal T^{\mathrm ord}_{w,b}$
be the sets of all pairs $(t,l)\in T^{\mathrm ord}$ such that
$t\in\mathcal T_P$ and $t\in\mathcal T_{w,b}$ respectively.
Given $s\in\mathcal T^{\mathrm ord}$ and $t\in\mathcal T$ we say
$s\in t$ if $s = (t,l)$ for some numeration $l$. This means that we
identify $t$ with the set of all pairs of the form $(t,l)$.

\nxpoint \label{s:counting}
(Example: counting trees in $\mathcal T^{\mathrm ord}$)
Let $s_w$ be the cardinality of $\mathcal T^{\mathrm ord}_{w,0}$, that is
the number of distinct numerated planar rooted trees with descending
numeration and only white nodes. We suggest reader to check that
$s_1 = s_2 = 1$, $s_3 = 3$ an $s_5 = 15$.
It is easy to derive a recursion for $s_w$. The trees in
$\mathcal T^{\mathrm ord}_{w+1,0}$ have a root node with at most $w$
numerated branches which are themselves planar rooted trees with labels.
The exact labelling is determined by first choosing
the set of labels of each branch, and then choosing a descending numeration
on the labels within each branch. For the whole process $w$ labels are
available, regarding that the root branch is mandatory labelled with $1$.
Thus we obtain the recursion
$$
s_{w+1} = \sum_{k = 1}^w \sum_{w_1+w_2+\ldots + w_k = w}
\frac{w!}{w_1!w_2!\cdots w_k!} s_{w_1} s_{w_2} \cdots s_{w_k},
\,\,\,\,\,\,w\geq 1.
$$
The solution of this recursion is
$s_w = (2w-3)!! = 1\cdot 3\cdot 5\cdots (2w-3)$.

Cardinality of $\mathcal T^{\mathrm ord}_{b,w}$ may be determined
similarly: for $w\geq 0$,
$$
s_{w+1,b} =  \sum_{k = 1}^w
\sum_{
\begin{array}{c}
w_1+\ldots + w_k = w\\
b_1 +\ldots + b_k = b
\end{array}}
\frac{w!}{w_1!w_2!\cdots w_k!} s_{w_1,b_1} s_{w_2,b_2} \cdots s_{w_k,b_k}.
$$

\nxpoint Suppose now $\gg$ and its basis $\hx_1,\ldots,\hx_n$ are fixed;
and hence the dual basis $\partial^1,\ldots,\partial^n$ and
the structure constants $C^i_{jk}$.

Given a tree $t\in\mathcal T^{\mathrm ord}_{w,b}$ and labels
$1\leq \mu,\alpha_1,\ldots,\alpha_w\leq n$, {\bf define}
$$\mathrm{ev}(t)^\mu_{\alpha_1,\ldots,\alpha_w}\in S(\gg^*)$$
as follows.
First replace the numeration labels $1,\ldots,w$ on white nodes
with $\alpha_1,\ldots,\alpha_w$. Then label arbitrarily
the inner lines by distinct new variables $\rho_1,\ldots, \rho_{w+b-1}$,
and attach a new {\it external} incoming line
to the root node and label it with label $\mu$.

To form an expression $\mathrm{ev}(t)^\mu_{\alpha_1,\ldots,\alpha_w}$
apply the {\it Feynman-like rules}: 
\begin{itemize}
\item To each {\bf white node}
with label $\alpha_k$, incoming node $\rho_l$ and
outgoing nodes $\rho_{v_1},\ldots,\rho_{v_s}$ assign
value $(-1)^s \frac{B_s}{s!}\sum_{k_1,\ldots,k_{s-1} = 1}^n
C^{k_1}_{\alpha_k\rho_{v_1}} C^{k_2}_{k_1 \rho_{v_2}}
\cdots C^{\rho_l}_{k_{s-1} \rho_{v_s}}$. \newline If $s = 0$ 
(the white node is a white leaf), the value is Kronecker delta
$\delta^{\rho_l}_\alpha$. 

\item To each {\bf black leaf} assign
$\partial^{\rho_l}\in\gg^*\subset S(\gg^*)$.

\item Multiply so assigned values of all nodes
and {\bf sum over all values from $1$ to $n$
of labels of all internal lines} 
$\rho_1,\ldots, \rho_{w+b-1}$.
\end{itemize}

\nxsubpoint Example:
\begin{equation}\label{eq:Feyrule}\xymatrix{
& \ar[d]^{\mu} & \\
& *++[o][F-]{\alpha_1} \ar[ld]_{\rho_1}\ar[rd]^{\rho_2} &
\\
{\bullet}&& *++[o][F-]{\alpha_2}
}
\frac{B_2}{2!} \sum_k (C^k_{\alpha_1 \rho_1}\partial^{\rho_1}) C^\mu_{k\rho_2}
\delta^{\rho_2}_{\alpha_2} = \frac{1}{12} \sum_k \CE^k_{\alpha_1} C^\mu_{k\alpha_2}
\end{equation}
\nxsubpoint
Clearly $\mathrm{ev}(t)^\mu_{\alpha_1,\ldots,\alpha_w}$ are components
of some tensor which will be of course denoted $\mathrm{ev}(t)
\in \gg\otimes T^n(\gg^*)$.
In this notation,
\begin{equation}\label{eq:commutatorevaltrees}
[\ldots[[\hpart^\mu,\hx_{\alpha_1}],\hx_{\alpha_2}],\ldots,
\hx_{\alpha_w}] =
\sum_{b = 0}^\infty\sum_{t\in \mathcal T^{\mathrm ord}_{w,b}}
\mathrm{ev}(t)^\mu_{\alpha_1,\ldots,\alpha_w}
\end{equation}

\nxpoint \label{s:fev}
For a tree $t\in\mathcal T^{\mathrm ord}_{w,b}$ one defines
its {\bf full evaluation}
$$
\mathrm{fev}(t)^\mu : =
\frac{1}{w!}\partial^{\alpha_1}\cdots\partial^{\alpha_w}\otimes
\mathrm{ev}(t)^\mu_{\alpha_1,\ldots,\alpha_w},
$$
and for $s\in\mathcal T$ one defines
$$
\mathrm{fev}(s)^\mu : = \sum_{t\in s, t\in\mathcal T^{\mathrm ord}}
\mathrm{fev}(t)^\mu .
$$
\nxpoint \label{s:symmselectionrule} {\bf (Basic selection rule)}
Suppose a tree $t\in\mathcal T$ has at least one white node $y$ such that
its most left child branch is a white leaf. Then for all $\mu$,
$$
\mathrm{fev}(t)^\mu = 0.
$$
{\it Proof.} Once the Feynman rules are applied the fact is rather obvious.
Namely, suppose that white node has $s$ child branches, its
label is $k$ and of its most left child branch is $l$ (then $l>k$).
Then the Feynman rules for
$\mathrm{ev}(t)^\mu_{\alpha_1,\ldots,\alpha_k,\ldots,\alpha_l,
\ldots,\alpha_w}$ assign to the white node $y$
the factor $(-1)^s \frac{B_s}{s!}C^*_{\alpha_1\rho_1} C^*_{*\rho_2}\cdots
C^{\rho_0}_{*\rho_s}$ if the incoming line to $y$ is labelled by $\rho_0$
and outgoing from left to right by $\rho_1,\ldots,\rho_s$. The white
leaf contributes by a factor $\delta^{\rho_1}_{\alpha_2}$. Thus we get
a factor of the type
$C^*_{\alpha_k\rho_1}\delta^{\rho_1}_{\alpha_l} = C^*_{\alpha_k\alpha_l}$
which is antisymmetric in lower indices.
To obtain $\mathrm{fev}(t)^\mu$ contract
$1\otimes\mathrm{ev}(t)^\mu_{\alpha_1,\ldots,\alpha_k,\ldots,\alpha_l,\ldots,
\alpha_w}$ with the symmetric tensor $\frac{1}{w!}\partial^{\alpha_1}\cdots
\cdots\partial^{\alpha_w}\otimes 1$ what vanishes by symmetry reasons. Q.E.D.
\vskip .02in
\nxsubpoint Notice that this selection rule holds for $\mathrm{fev}$ but not
for $\mathrm{ev}$ (the latter does not involve symmetrization). The
subset of trees which are not excluded in calculation of $\mathrm{fev}$
by the basic selection rules are called
$\mathrm{(fev)}$-{\bf contributing trees} and the
correspoding subclasses are distingushed with supersctipt $c$,
e.g. $\mathcal T^c_{w,b}\subset \mathcal T_{w,b}$.

\vskip .04in

By similar symmetry reasons, the following result holds:

\nxpoint \label{s:symmonvacuum} {\bf Lemma.} {\it
Let $\hx_1,\ldots,\hx_n$ be a basis of $\gg$.
If $\xi : S(\gg)\to U(\gg)$ is the coexponential map, then for $w\geq 2$,
$$
\sum_{\sigma\in\Sigma(w)}
[\ldots[[\hpart^\mu,\hx_{\sigma \alpha_1}],\hx_{\sigma \alpha_2}],\ldots,
\hx_{\sigma\alpha_w}]
(1) = 0,
$$
where on the left hand side
the evaluation at unit element (``vacuum'') is in the sense of
the action of the Weyl algebra on the usual
symmetric algebra $S(\gg)$.
}

The evaluation at vacuum simply kills all the strictly positive powers
of $\partial$-s, hence only the terms coming from trees in
$\mathcal T^{\mathrm ord}_{w,0}$ survive. Thus the lemma may be restated as
$$
\sum_{\sigma\in\Sigma(k)} \sum_{t\in\mathcal T^{\mathrm ord}_{w,0}}
\mathrm{ev}(t)^\mu_{\sigma\alpha_1\cdots\sigma\alpha_w} = 0.
$$
The proof in the latter form is obvious:
applying the Feynman rules to a graph with
$w$ nodes and $w-1$ internal lines
produces a tensor
which is proportional to some
contracted product of $w-1$ copies of the structure constants tensor $C$,
$w-1$ contractions, $w$ lower external labels
and one upper external label $\mu$. In particular at least one pair of
labels $\alpha_i,\alpha_j$ will be attached as lower labels
of the same $C$-tensor. By the antisymmetry in subscripts of $C$, after
symmetrization of $\alpha_1,\ldots,\alpha_w$ we obtain zero.

\nxpoint \label{s:symmderanclassical} {\bf Corollary.} {\it
In the symmetric ordering (if $\xi$ is the coexponential map),
the formula for the derivatives of $(\hat{a})^p = (a^\beta\hx_\beta)^p$
is of the classical (undeformed) shape, i.e.
$$
\frac{1}{s!}\hpart^{\alpha_1}\hpart^{\alpha_2}\cdots\hpart^{\alpha_s}
(\hat{a}^p)
= {p\choose s}
a^{\alpha_1} a^{\alpha_2}\ldots a^{\alpha_s} \hat{a}^{p-s}, \,\,\,\,p\geq s.
$$
}
This follows by an induction on $k$; the induction step involves applying the
case $k =1$.
For $k =1$, the formula follows from~(\ref{eq:dxn}) for $\hat{f} = 1$
after noticing that $a^{\alpha_1}a^{\alpha_2}\cdots a^{\alpha_k}$
in~(\ref{eq:dxn}) is symmetric
under permutations of $\alpha_1,\ldots,\alpha_k$, hence by
\refpoint{s:symmonvacuum} the only term which survives is the top degree
term which is of classical shape.

\nxpoint Up to the fourth order in total derivative,
or equivalently, third order in $C$-s one gets the following
$$\begin{array}{lcl}
\Delta \hpart^\mu &=& 1 \otimes \hpart^\mu + \hpart^\mu \otimes 1
+ \frac{1}{2} C^\mu_{\alpha\beta} \hpart^\alpha\otimes \hpart^\beta
+ \frac{1}{12} C^\star_{\alpha\beta} C^\mu_{\star\gamma}
(\hpart^\alpha \otimes \hpart^\beta\hpart^\gamma +
\hpart^\beta \hpart^\gamma \otimes
\hpart^\alpha) \\
&& \,\,\,-\frac{1}{24}C^*_{\alpha\beta}C^*_{*\gamma}C^\mu_{*\delta}
\hpart^\alpha\hpart^\gamma \otimes\hpart^\beta\hpart^\delta + O (C^4)
\end{array}$$
where we sum on pairs of repeated indices (including $*$,
where on two consecutive ones).

\nxpoint {\bf Theorem.} \label{s:coprodbrabra} {\it
If $\xi$ is the coexponential map, the coproduct is given by
$$
\Delta\hpart^{\mu} = 1\otimes\hpart^\mu +
\hpart^\alpha\otimes[\hpart^\mu,\hx_\alpha]
+\frac{1}{2} \hpart^\alpha\hpart^\beta\otimes
[[\hpart^\mu,\hx_\alpha],\hx_\beta] +\ldots
$$
or, in symbolic form,
\begin{equation}\label{eq:expeuler}
\Delta\hpart^{\mu} = \exp(\hpart^\alpha\otimes \mathrm{ad}\,(-\hx_\alpha))
(1\otimes\hpart^\mu)
\end{equation}
and in the tree expansion form, using the notation from~\refpoint{s:fev},
\begin{equation}\label{eq:Deltafev}
\Delta\hpart^{\mu} = \sum_{t\in\mathcal T^{\mathrm ord}} \mathrm{fev}(t)^\mu.
\end{equation}
}
Of course, each $\ad(-\hx_\alpha)$ in~(\ref{eq:expeuler})
has to be applied to $\hpart^\mu$ before applying
the whole expression on the elements in $S(\gg)\otimes S(\gg)$
(for the Leibniz rule for the star product)
or on the elements in $U(\gg)\otimes U(\gg)$
(for the Leibniz rule for the usual noncommutative product).

{\it Proof.} It is well known that the expressions of the form
$(\hat{a})^p$ where $\hat{a} = \sum_\alpha a^\alpha \hx_\alpha$ with
varying $a = (a^\alpha)$ span $U(\gg)$. Thus it is sufficient to show
that for all $a$, all $\hat{f}\in U(\gg)$ and all $p$ the
twisted Leibniz rule
$$
\hpart^\mu(\hat{a}^p \hat{f}) = \sum_{w=0}^p \frac{1}{w!}
\sum_{\alpha_1,\ldots,\alpha_w}
\hpart^{\alpha_1}\cdots\hpart^{\alpha_w}(\hat{a}^p)
[[\ldots[\hpart^\mu,\hx_{\alpha_1}],\ldots],\hx_{\alpha_w}](\hat{f}).
$$
holds. This follows by comparing the Corollary
\refpoint{s:symmderanclassical}
which holds for symmetric ordering only
with the formula~(\ref{eq:dxn}) which holds for general ordering.

\nxpoint
Let $\partial^{a b c} = \partial^a\partial^b\partial^c $ and so on.
Recall $\phi^\mu_\nu = \phi^\mu_\nu(\partial) = [\hpart^\mu,\hx_\nu]$.

{\bf Corollary.} {\it In symmetric ordering, 
for any $\hat{f},\hat{g}$ in $U(\gg)$,}
$$
\phi^\mu_\nu (\partial)(\hat{f}\hat{g}) = \sum_{N = 1}^\infty
\frac{1}{N!}\sum_{i_1,\ldots,i_N}\sum_{k= 1}^N
\partial^{i_1\cdots i_{k-1} i_{k+1}\cdots i_N}
\phi^{i_k}_\nu(\partial)
(\hat{f})\cdot
[[\ldots [\partial^\mu,\hat{x}_{i_1}], \ldots, \hat{x}_{i_N}]
(\hat{g})
$$

Notice that the last sum is from 1, not 0.
Summation over repeated indices understood.
This formula is equivalent to giving
the formula deformed coproduct for the argument
$\Delta([\hpart^\mu,\hx_\nu]) = \Delta(\phi^\mu_\nu)$.
For the proof, calculate
$\hpart^\mu((\hx_\nu \hat{f})\hat{g})$ using the twisted Leibniz rule
from the theorem \refpoint{s:coprodbrabra}, and subtract
similarly $\hx_\nu \hpart^\mu (\hat{f}\hat{g})$ and group the terms
and commutators appropriately.

\nxpoint \label{s:s1ptp1symmetry}
Let $\tau : S(\gg^*)\hat\otimes S(\gg^*)\to  S(\gg^*)\hat\otimes S(\gg^*)$
be the standard flip interchanging the tensor factors (in the completed
tensor product).

{\bf Theorem.} {\it
Let $s_{1,p}$ be the unique tree in $\mathcal T^{\mathrm ord}_{1,p}$.
Then for all $\mu$,
\begin{equation}
\tau(\mathrm{fev}(s_{1,p})^\mu) =
(-1)^{p+1}\sum_{(t,l)\in\mathcal T^{\mathrm ord}_{p,1}}\mathrm{fev}(t)^\mu
\end{equation}
or more explicitly
\begin{equation}\label{eq:symmetry1pexplicit}
\mathrm{ev}(s_{1,p})^\mu_\beta \otimes \partial^\beta
= (-1)^{p+1}\sum_{t\in\mathcal T^{\mathrm ord}_{p,1}}
\frac{1}{p!} \partial^{\alpha_1}\cdots\partial^{\alpha_p}
\otimes \mathrm{ev}(t)^\mu_{\alpha_1,\ldots,\alpha_p}
\end{equation}
}
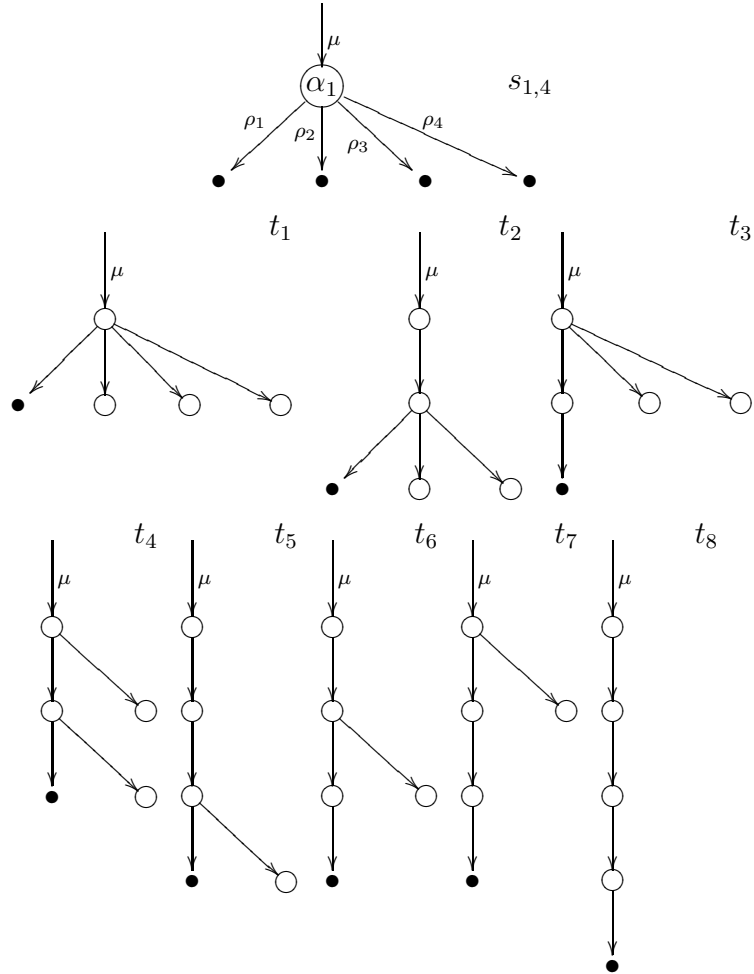
$$\xymatrix{
& \ar[d]^{\mu} & & \\
& *++[o][F-]{\alpha_1} \ar[ld]_{\rho_1}\ar[d]_{\rho_2}\ar[rd]_{\rho_3}
\ar[rrd]^{\rho_4} & & s_{1,4}
\\
\bullet&\bullet&\bullet&\bullet
}
$$
$$
\xymatrix{
& \ar[d]^{\mu} & & t_1\\
& *++[o][F-]{} \ar[ld]\ar[d]\ar[rd]
\ar[rrd]& &
\\
\bullet&*++[o][F-]{}&*++[o][F-]{}&*++[o][F-]{}
}\,\,\,\,
\xymatrix{
& \ar[d]^{\mu} &  t_2\\
& *++[o][F-]{}\ar[d]&\\
& *++[o][F-]{} \ar[ld]\ar[d]
\ar[rd]&
\\
\bullet&*++[o][F-]{}&*++[o][F-]{}
}\,\,\,\,
\xymatrix{
\ar[d]^{\mu} & & t_3\\
*++[o][F-]{}\ar[d]\ar[rd]\ar[rrd]&\\
*++[o][F-]{}
\ar[d]&*++[o][F-]{}&*++[o][F-]{}
\\
\bullet&&
}
$$
$$
\xymatrix{
\ar[d]^{\mu} &  t_4\\
*++[o][F-]{}\ar[d]\ar[rd]&\\
*++[o][F-]{}\ar[d]
\ar[rd]& *++[o][F-]{}
\\
\bullet&*++[o][F-]{}
}\,\,\,
\xymatrix{
\ar[d]^{\mu} &t_5\\
*++[o][F-]{}\ar[d]\\
*++[o][F-]{}\ar[d]\\
*++[o][F-]{}\ar[d]\ar[rd]
\\
\bullet & *++[o][F-]{}
}\,\,\,
\xymatrix{
\ar[d]^{\mu} &t_6\\
*++[o][F-]{}\ar[d]\\
*++[o][F-]{}\ar[d]\ar[rd]\\
*++[o][F-]{}\ar[d]& *++[o][F-]{}
\\
\bullet &
}\,\,\,
\xymatrix{
\ar[d]^{\mu} &t_7\\
*++[o][F-]{}\ar[d]\ar[rd]\\
*++[o][F-]{}\ar[d]& *++[o][F-]{}\\
*++[o][F-]{}\ar[d]
\\
\bullet
}\,\,\,
\xymatrix{
\ar[d]^{\mu} &t_8\\
*++[o][F-]{}\ar[d]\\
*++[o][F-]{}\ar[d]\\
*++[o][F-]{}\ar[d]\\
*++[o][F-]{}\ar[d]
\\
\bullet
}
$$
The diagrams above show $s_{1,4}$ and the 8 diagrams
$t_1,\ldots,t_8\in\mathcal T_{4,1}^c$.

{\it Proof.} For $p = 1$ the assertion is a tautology. Let us prove
the assertion for $p>1$.

By the Feynman rules,
the LHS of~(\ref{eq:symmetry1pexplicit}) equals
$$
(-1)^p \frac{B_p}{p!} \sum_{k_1,\ldots,k_{p-1}}
C^{k_1}_{\beta\rho_1} C^{k_2}_{k_1\rho_2}
\cdots C^\mu_{k_{p-1}\rho_p}
\partial^{\rho_1}\partial^{\rho_2}\cdots\partial^{\rho_p}
\otimes\partial^\beta
= (-1)^p \frac{B_p}{p!}(\CE^p)^\mu_\beta \otimes \partial^\beta.
$$
Therefore it is sufficient and we will show by induction that
$$
\sum_{t\in\mathcal T^{\mathrm ord}_{p,1}}
\frac{1}{p!} \partial^{\alpha_1}\cdots\partial^{\alpha_p}
\otimes \mathrm{ev}(t)^\mu_{\alpha_1,\ldots,\alpha_p}
= (-1)^{p+1} \frac{B_p}{p!}(\CE^p)^\mu_\beta \otimes \partial^\beta.
$$
Bournulli numbers and hence this expression are zero for odd $p>1$
and nonzero for even $p>1$.

By the basic selection rule \refpoint{s:symmselectionrule}, the only
trees $t\in \mathcal T_{p,1}$ (no labelling) which may
give a nonzero contribution are those who have no leftmost white leafs,
and regarding that there is only one black node in our case, only
one white node may have a leftmost leaf (which is black).
That means that every contributing tree in $\mathcal T_{p,1}^c$ is
composed as follows: start with a vertical chain made of $r+1\leq p$
white nodes ending with a black node on the bottom and on
this white chain there are attached $(p-r-1)\geq 0$
right-hand side leafs (to some among the white nodes of the
vertical chain), but no branches of length $\geq 2$ are attached.

Notice that each $t\in \mathcal T_{p,1}^c$ for $p>1$  may be also composed
alternatively starting with the top white node,
attaching the left-most branch $t'\in\mathcal T^c_{r,1}$
and $p-r-1$ leafs, $r\geq 0$. We group the trees by the number $0\leq r<p$.
Let us now consider the ordered trees $t\in \mathcal T_{p,1}^{c,\mathrm{ord}}$.
To the top node we must assign label $1$, then we may choose any $r$
remaining numbers $\beta_1,\ldots,\beta_r$
to distribute them within $t'$ branch according to the
usual ordering rules within $t'$ and distribute the remaning $p-r-1$ labels
$\gamma_1,\ldots,\gamma_{p-r-1}$
to the white leafs in any order. Other way around,
given $t$ with labels, if $t'$ as a branch of $t$,
then its labels are renumerated as $1$ to $r$
in the same order. For example,
labels $2,5,7,8,3$ of white nodes in $t'$ as a branch
will be replaced by the position labels $1,3,4,5,2$
in $t'$ as an independent tree.
Thus for a given ordering
$$
\mathrm{ev}^\mu_{1,\ldots,r} (t)
= (-1)^p\frac{B_p}{p!} \sum_{\rho,k_1,\ldots,k_{p-r-1}}
C^{k_1}_{1,\rho} C^{k_2}_{k_1\gamma_1}\cdots
C^\mu_{k_{p-r-1}\gamma_{p-r-1}} \mathrm{ev}^\rho_{\beta_1,\ldots,\beta_r}(t')
$$
(of course each $i$ has to be replaced by $\alpha_i$). Now we need to
count all ordering and combine into $\mathrm{fev}$. The ordering
constraints described above give some combinatorial factors, as well
as $1/n!$ in the definition of $\mathrm{fev}$. We obtain
$$
\sum_{t\in\mathcal T^{c}_{1,p}}
\mathrm{fev}(t)^\mu =
\frac{1}{p!} \sum_{r}(-1)^{p-r+1} \frac{B_{p-r}}{(p-r)!} {p-1\choose r} r!
\,\mathrm{fev}(t')^\rho(p-r-1)!((\CE^{p-r})^\mu_\rho\otimes  1).
$$
Notice here an additional sign from the first $C$-factor
(by antisymmetry of lower indices):
$C^*_{\alpha_1\rho}\partial^{\alpha_1} = -\CE^*_\rho$.

By the induction hypothesis, $\mathrm{fev}(t')^\rho =
(-1)^r\frac{B_r}{r!}(\CE^r)^\rho_\beta\otimes \partial^\beta$, hence,
$$
\sum_{t\in\mathcal T^{c}_{1,p}}
\mathrm{fev}(t)^\mu =
\frac{1}{p}(-1)^{p} \sum_{r} \frac{B_{p-r}}{(p-r)!} \frac{B_r}{r!}
((\CE^p)^\mu_\beta\otimes  \partial^\beta)
$$
Regarding that, for $p>1$,
$B_{p-r}$ and $B_r$ on the right are simultaneously nonzero
if and only if $r$ and $p-r$ are both even,
the proof finishes by applying the well known identity for Bernoulli numbers
$$
\sum_{s=1}^{l} \frac{B_{2s}}{(2s)!}\frac{B_{2l-2s}}{(2l-2s)!}
= \frac{-B_{2l}}{(2l-1)!} +\frac{1}{4}\delta_{l,1},\,\,\,\,\,\,\,\,\,l>0.
$$

\nxpoint
In these terms we state the following conjecture on the star product

In our notation we will often not distinguish any more $\partial$ from
$\hat\partial$; with the convention that when we write $[\partial,\hx]$
where $\hx\in U(\gg)$ we mean $\hat\partial$;
as well as when we apply $\hat\partial(\hat{f})$ with $\hat{f}\in U(\gg)$;
however when we apply $\partial(f)$ with $f\in S(\gg)$ we mean the usual
(undeformed) Fock representation. In any case $\Delta$ is
deformed and $\Delta_0$ undeformed coproduct: $\Delta_0(\partial^\mu)
= 1\otimes\partial^\mu + \partial^\mu\otimes 1$.

\nxpoint \label{s:conj} {\bf General conjecture.} (for all $\bm\phi$)
\begin{equation}\label{eq:conj2}
f \star g = \sum_{i_1,i_2,\ldots,i_n \geq 0}^{}
\frac{x^{i_1}_{1} x^{i_2}_{2}\cdots x^{i_n}_{n}}{i_1!\cdots i_n!}\,
m \left(\left(\prod_{l = 1}^n (\Delta - \Delta_0)((\partial^{l})^{i_l})\right)
(f\otimes g)\right),
\end{equation}
where $f,g \in S(\gg)$ and $m$ is the commutative
multiplication of polynomials
$S(\gg)\otimes S(\gg)\to S(\gg)$. Notice that for any concrete $f$ and $g$,
the summation on the right has only finitely many nonzero terms.
This formula is proved in some special cases~(\cite{MS}) and in this article
for general $\gg$ and symmetric ordering. For general $\bm\phi$, if
$f$ is a first order monomial and $g$ arbitrary, this formula boils down
to our main formula of article~\cite{ldWeyl}.

Formula~(\ref{eq:conj2}) can be expressed via normal ordered exponential
:exp(): (here $x$-s to the left, $\partial$-s to the right)
$$
f\star g = m : \exp(x_\alpha (\Delta -\Delta_0)(\partial^\alpha)): (f\otimes g)
$$
and $m$ is the usual product. Notcie that 
$: \exp(x_\alpha (\Delta -\Delta_0)(\partial^\alpha)):$ is not 
an element of the tensor product $H\otimes H$ where $H$ is the algebra of
formal vector fields;
namely the position of $x$-variables is to the left from the 
$\partial$-s, but the tensor factor is not chosen, and it does not
matter as we use $m$ after application of the derivatives to $f\otimes g$.
But we believe there is a correct alternative form where the positions
of $x$-s in tensor factor is chosen and additional good properties
(yielding a Drinfeld twist) are satisfied (cf. ~\refpoint{s:cases}).

\nxpoint In articles~\cite{covKappadef,MS} for a particular Lie algebra,
the case of ``kappa-deformed Euclidean space'' in dimension $n$ 
for which the commutation relations are of the type
$[\hx_\mu,\hx_\nu] = i(a_\mu \hx_\nu - a_\nu\hx_\mu)$
for some vector $\bm a = (a_1,\ldots, a_m)$,
the conjecture has been verified for general $\phi$.

\nxpoint \label{s:mainthm} {\bf Main theorem.} {\it For symmetric ordering
the conjecture \refpoint{s:conj} holds for all $\gg$.
}

In fact we can prove the conjecture 
in more general case, for those $\bm\phi$
which are obtained using certain procedure of twisting basis
by a wide class inner automorphisms of semicompleted Weyl algebra.

\section {Some facts on Hausdorff series}

\nxpoint {\bf (The recursive form of Hausdorff series)}
Given $X,Y\in\gg$ where $\gg$ is finite-dimensional
with a norm inducing the standard topology.
The series $H(X,Y)$ is uniquely defined by
$$
\exp(X)\exp(Y) = \exp(H(X,Y))
$$
and it converges in such norm.
Then $H(X,Y) = \sum_{N=0}^\infty H_N(X,Y)$ where
{\bf ``Dynkin's Lie polynomials''} $H_N = H_N(X,Y)$ are defined
recursively by $H_1 = X+Y$ and
$$
(N+1)H_{N+1} = \frac{1}{2}[X-Y,H_N] +
\sum_{r = 0}^{\lfloor N/2 -1\rfloor}
\frac{B_{2r}}{(2r)!}\sum_{
s\nodo{
\begin{array}{c} s_1 +\ldots + s_{2r} = N\\
s_i > 0, 1\leq i \leq 2r
\end{array}}
}
[H_{s_1},[H_{s_2},[ \ldots, [H_{s_{2r}},X+Y]\ldots]]]
$$
where the sum over $s$ is the sum over all
$2r$-tuples $s = (s_1,\ldots,s_{2r})$
of strictly positive integers
whose sum $s_1 +\ldots+s_{2r} = N$. This identity is
well-known and we do not reprove it here.

\nxpoint {\bf (Linear parts in either $X$ or $Y$)} 
The linear part in $X$ of the Hausdorff series is
$H_{1,\star}(X,Y) = \sum_{N=0}^\infty (-1)^N \frac{B_N}{N!} [Y,[\ldots,[Y,X]]]$
where $N$ is the degree of $Y$ in the Lie polynomial involved.
Similarly, the linear part in $Y$ is
$H_{\star,1}(X,Y) =  \sum_{N=0}^\infty \frac{B_N}{N!} [X,[\ldots,[X,Y]]]$
where $N$ is the degree of $Y$ in the Lie polynomial involved.

\nxpoint {\bf (Symmetries of Hausdorff series)} Identity
$e^X e^Y = (e^{-Y} e^{-X})^{-1}$ implies $H(-Y,-X) = -H(X,Y)$.
Dynkin's polynomials are of fixed total degree, hence the change
$(X,Y)\mapsto (-Y,-X)$ does not mix them and $H_P(-Y,-X) = -H_P(X,Y)$
for all $P>0$. We refine the degree grading on a free Lie algebra
on two generators by a bigrading which induces a decomposition
$H_P(X,Y) = \sum_{w+b=P} H_{w,b}(X,Y)$ where $H_{w,b}$
is the sum of all Lie polynomials in $H_P(X,Y)$
of degree $w$ in $X$ and degree $b$ in $Y$. Clearly, knowing
$H_P$ determines $H_{w,b}$ for all $w,b$ with $w+b = P$.

\nxpoint {\bf Proposition.}
{\it The following $w$-recursion and $b$-recursion hold
$$
(w+1)H_{w+1,b} = \frac{1}{2}[X,H_{w,b}] +
\sum_{r = 0}^{\lfloor w/2 -1\rfloor}
\frac{B_{2r}}{(2r)!}\sum_{
w_i, b_i}
[H_{w_1,b_1},[\ldots, [H_{w_{2r},b_{2r}},X]\ldots]]
$$
$$
bH_{w+1,b} = -\frac{1}{2}[Y,H_{w,b}] +
\sum_{r = 0}^{\lfloor b/2 -1\rfloor}
\frac{B_{2r}}{(2r)!}\sum_{
w_i, b_i}
[H_{w_1,b_1},[\ldots, [H_{w_{2r},b_{2r}},Y]\ldots]]
$$
where in the sum on the RHS $\sum_i w_i = w$ and $\sum_i b_i = b$
for the $w$-recursion and $\sum_i w_i = w+1$ and $\sum_i b_i = b-1$
for the $b$-recursion.
}

{\it Proof.} For the purpose of the proof we introduce two
new sets of Lie polynomials. The first set will have members
$H^W_{w,b}$ and the latter $H^B_{w,b}$ where $w\geq 0, b\geq 0, w+b>0$.
For $w =0$ we set $H^W_{w,b} = H(w,b)$ what is $0$ unless $b = 1$ when
$H^W_{0,1}= X$; similarly for $b = 0$
we set $H^B(w,b)=H(w,b)$. Also set $H^W_{1,0} = Y$  and $H^B_{0,1} = X$,
regarding that $(0,0)$ point is undefined.
By definition, $w$-recursion is used to define
$H^W_{w,b}$ at all other pairs $(w,b)$ and similarly the $b$-recursion
is used to define $H^B_{w,b}$. E.g. for $w$-recursion we first
use the recursion at the line $b = 0$, increasing from $w =1$ on,
then at the line $b =1$, increasing from $w =1$, and so on. Clearly each
recursion relation is used exactly once to determine one new value and
all instances of relations are used. Notice that on the line $b = 0$,
the $w+1 = w+b+1 = P+1$, hence the $w$-recursion gives the same values on
this line as the standard recursion for $H_{w,b}$. In that manner we
notice that the initial values (line $b = 0$ and $(0,1)$) given to $H^B$
agree with the value of $H^W$ and $H$ obtained by $w$-recursion and the
standard recursion. The initial values hence also satisfy the symmetries
$H_{w,b}(X,Y) = -H_{b,w}(-Y,-X)$ in both cases.
We want to prove that the values within the quadrant agree as well,
not only the conditions on the boundary.
But, the $b$-recursion may be obtained from $w$-recursion
also by the same symmetry operation! Regarding that the symmetry holds
for initial values and also for the recursion, than this is true for
each pair of new points to which the two recursions assign the values.
Conclusion: $H^B = H^W$. Therefore we can now safely combine two recursions
without being afraid of nonconsistency. But adding up the $w$-recursion and
$b$-recursion we clearly get the standard recursion. Regarding that the
initial value $w+b = P = 1$ for standard recursion is checked and that the
standard recursion is the consequence, and also that the values $H_P$
determine $H_{w,b}$, we conclude $H = H^B = H^W$.

\nxpoint {\bf (Recursive formula for $D = D(k,q)$)}
Let $\hx_1,\ldots,\hx_n$ be a basis of $\gg$, $i = \sqrt{-1}$,
$X = i k^a \hx_a$, $Y = i q^a \hx_a$ and $H(X,Y) = iD^a(k,q)\hx_a$,
where $k = (k^1,\ldots,k^n)$, $q = (q^1,\ldots,q^n)$; let also
$D = D(k,q) = (D^1(k,q),\ldots,D^n(k,q))$.
Then $D^\mu(k,q) = \sum_{N=0}^\infty D^\mu_n(k,q)$
where $D_1^\mu(k,q) = k^\mu + q^\mu$
and the recursion
\begin{equation}
(N+1) D^\mu_{N+1} = \frac{1}{2}(k^a-q^a)(E_N)^\mu_a +
\sum_{r = 1}^{\lfloor N/2 -1\rfloor} \frac{B_{2r}}{(2r)!}
\sum_s (k^a + q^a)(E_{s_1}\cdots E_{s_{2r}})^\mu_a
\end{equation}
holds where
$$
(E_P)^\mu_\nu := \sum_\sigma iC^\mu_{\nu\sigma} D^\sigma_P, \,\,\,\,\,
P\geq 1,
$$
are the components of a matrix $E_P$, and the product of matrices on
the right is via the convention that the superscript is the {\it row} index.
The sum over $a$ on the right is understood and the sum over $s$
is again over $2r$-tuples of positive integers adding up to $N$.

\section{Fourier notation and using exponentials}

\nxpoint \label{s:Fourier} 
If some linear isomorphism $S(\gg)\to U(\gg)$
preserves the degree filtration, then it clearly extends by continuity
to a linear map among the corresponding
completions $\widehat{S(\gg)}\to\widehat{U(\gg)}$.
If the isomorphism is a coalgebra map, then the extension respects the
completed coproducts $\Delta : \hat{H}\to \hat{H}\hat\otimes\hat{H}$
($H = S(\gg)$ or $U(\gg)$). Thus, it makes sense to consider the behaviour of
exponential series (as a formal series) under coalgebra isomorphism
$\xi$ as above. It is also useful to extend the field by $\sqrt{-1}$ if
it is not present and consider formal series of the type
$\exp(ik^\alpha x_\alpha)$. If the field is $\mathbb C$ then such
series are specially important because of  Fourier integral methods.
However, Fourier integral is defined only for some formal series, so
the formulas, though useful for other spaces of functions (one can
extend our coproducts etc. to various functional spaces, but
we will avoid this here) the formulas involving Fourier integrals in this
paper will be understood just in the following sense:
every abstract series involved is a finite sum of formal power series of
the form $\exp(ia^\alpha x_\alpha)$. The linear space of such
such finite sums (of exponentials), $S_e(\gg) \subset \widehat{S(\gg)}$
is dense in the space of all formal power series. Thus
if we prove that some identity between functionals continuous with respect to
the power series filtration, holds when restricted to this space, 
the identity holds in general.  
Even when the identity is proved
for finite sums of exponentials we heuristically write
integrals, instead of sums.
The imaginary unit is just for suggestiveness of applications in physics,
one can correct the $\sqrt{-1}$ factors and prove the formulas
just for the sums of functions of the form $\exp(ia^\alpha x_\alpha)$
but we will not spend time on these nicetess.

\nxpoint
Coalgebra isomorphisms $\xi : S(\gg)\to U(\gg)$
which are identity on $\genfd\oplus\gg$, and which are 
extended to the completions have the property
\begin{equation}\label{eq:Kfcn}
\xi(\exp(ik^\alpha x_\alpha)) = \exp(iK(\vec{k})^\beta \hx_\beta)
\end{equation}
for some bijection $K : \genfd^n\to\genfd^n$. 
(Proof: All group like elements both in $\hat{S}(\gg)$ and
in $\hat{U}(\gg)$ are of such exponential form. $\xi$ is a bijection 
and preserves the group like elements because it is a coalgebra map.)
For example, if $K$ is the identity map, 
this is the case of symmetric ordering:
$\xi$ is the coexponential map (when considered defined on
$S(\gg)$ only). Furthermore, one can get a
very large class of other solutions
which satisfy~(\ref{eq:Kfcn})
using certain inner automorphisms of Weyl algebra.

Namely, let $S = \exp(x^\alpha R_\alpha + B)$
where $R_\alpha = R_\alpha(\partial), B= B(\partial)$ are
some formal series in variables
$\partial^1,\ldots,\partial^n \in \hat{A}_{n,\genfd}$.
Then define the formal power series
$$
y_\alpha := Sx_\alpha S^{-1},\,\,\,\,\,\,\,\,\,\,\,
\partial^\alpha_y  := S\partial^\alpha S^{-1}.
$$
They again satify canonical commutation relations:
$\partial^\alpha_y,y_\beta] = \delta^\alpha_\beta$
(this does not depend on the special form of $S$) and
$y_\alpha = x_\rho \psi^\rho_\alpha$ for some formal power series
in $\partial$-s
$\psi^\rho_\alpha = \psi^\rho_\alpha(\partial)$
(this follows by the special form of $S$). Moreover,
$ \partial^\alpha_y = d(\partial)$ is also a power series
in $\partial$-s only.

Now $x_\alpha \phi^\alpha_\beta = y_\rho \psi^\rho_\tau \phi^\tau_\beta$.
This way $\phi^\alpha_\beta$ is in the new basis replaced by
$\psi^\rho_\tau(\partial) \phi^\tau_\beta(\partial)$
what should be expressed in terms of
$\partial_y$-s (what is done computing the inverse transformation $S^{-1}$).
This way we get some $\Xi^\alpha_\beta = \Xi^\alpha_\beta (\partial_y)$
in place of $\phi^\alpha_\beta$.
This procedure can be accomplished in some special cases
(for $\kappa$-Minkowski space see~\cite{kappaind,MS}),
but also some general statements may be proved for $\Xi^\alpha_\beta$
obtained by this procedure, where the starting $\phi^\alpha_\beta$ 
corresponds to the symmetric ordering.

\section{The results leading to the proof of the main theorem.}

\nxpoint For the coexponential map $\xi$, the equality
$\xi(\exp(ik^a x_a)) = \exp(ik^a\hx_a)$ holds.
Therefore the star product $f\star g = \xi^{-1}(\xi(f)\cdot\xi(g))$
reduces to calculations with Hausdorff series. Namely
if $f(x) = \exp(ik^a x_a)$, $g(x) = \exp(iq^a x_a)$, then
$(f\star g)(x) = \exp(iD^a(k,q)x_a)$. For general $f$ and $g$,
it is convenient to expand $f$ and $g$ in Fourier components
(reasoning understood in the sense of~\refpoint{s:Fourier})
$f(x) = \int \frac{d^n k}{(2\pi)^n} (Ff)(k) \exp(ik^a x_a)$ and,
by bilinearity, we obtain
$$
(f\star g)(x) = \int \frac{d^n k}{(2\pi)^n} \int \frac{d^n k}{(2\pi)^n}
(Ff)(k) (Fg)(q) \exp(iD^a(k,q)x),
$$
or alternatively,
$$
(f\star g)(x) = m \exp(iz_a(D^a(-i\partial\otimes 1,-i\otimes\partial)+
i\partial^a\otimes 1 +i\otimes\partial^a)) (f\otimes g)(x)|_{z_a = x_a}
$$
where $\partial = (\partial^1,\ldots,\partial^n)$. Now notice that
$D^a_1(-i\partial\otimes 1,-i\otimes\partial) =
-i\partial^a\otimes 1 -i\otimes\partial^a$, hence
$$
(f\star g)(x) =
m \exp(iz_a(D^a-D^a_1)(-i\partial\otimes 1,-i\otimes\partial))
(f\otimes g)(x)|_{z_a = x_a}
$$
Notice that $iD^a_1(-i\partial\otimes 1,-i\otimes\partial)
= \Delta_0(\partial^a)$.
In fact, using the filtration by the total degree, we now see that the
{\bf main theorem \refpoint{s:mainthm} is equivalent to}

\nxpoint \label{s:t2} {\bf Theorem.} {\it
Let $\Delta_P(\partial^a)$ be the summand in $\Delta(\partial^a)$
consisting of terms of total homogeneity $P\geq 1$. Then for every $P\geq 1$,
$$
iD^a_P(-i\partial\otimes 1,-i\otimes\partial) =
\Delta_P(\partial^a)
$$
}
The theorem will be proved by induction on $P$. In other words, we have to
prove the corresponding recursion for $\Delta_P$.
We use two tools: 1. Fourier transform (this is only heuristic term here,
strictly speaking we use the denseness of the linear span of all
exponential series $\exp(a^\alpha x_\alpha)$ in $\widehat{S(\gg)}$
and do not require the existence of the imaginary unit,
as explained in~\refpoint{s:Fourier})
and 2. the combinatorics of the trees
whose Feynman rule contribution is involved here. 
Namely expression~(\ref{eq:Deltafev})
can be recursively computed after being filtered by the bidegree
$$
\Delta\hpart^{\mu} = \sum_{t\in\mathcal T^{\mathrm ord}} \mathrm{fev}(t)^\mu
= \sum_{b+w = 1}^\infty \sum_{t\in\mathcal T^{\mathrm ord}_{b,w}}
\mathrm{fev}(t)^\mu 
= \sum_{b+w = 1}^\infty \Delta_{b,w}\hpart^\mu.
$$
After evaluating, $b$ and $w$ correspond to the power of $\partial$-s
in the left and right tensor product factor respectively.
In other words, exactly 
Every degree in homogeneity corresponds to a node (white nodes
for right tensor product factor and black leaves for left factor)
as it is easily seen from the expression for $\mathrm{fev}$ and
Feynman rules for $\mathrm{ev}$.
The new node in induction procedure can always be assumed to be the top node, 
and, in particular white. Then one uses the two crucial lemmas
which use the expansions encoded in our Feynman-rule calculus:

\nxpoint {\bf Lemma.} {\it w-recursion formula holds 
(in Fourier transformed form) for calculating 
$\Delta\hpart^{\mu}$, where increasing w-degree by 1 corresponds to
one white node added and this w-recursion formula
is the same as for Hausdorff series in Fourier trasnformed form.
}

The proof is obtained using our Feynman rules and 
accounting for correct combinatorial factors,
in the same way as the counting of trees in~\refpoint{s:counting}, 
but with weights. We leave this to the reader.

\nxpoint {\bf Lemma.} {\it The initial conditions for w-recursion are the
same as for the w-recursion of the Hausdorff series}.

This lemma follows from \refpoint{s:s1ptp1symmetry}. 

Therefore the theorem \refpoint{s:t2} follows and
hence the main theorem \refpoint{s:mainthm}. 

\section{Special cases and other results}

\nxpoint \label{s:cases} {\bf (Classical cases and twists)}
The classical case of Moyal noncommutative space, where 
the deformation is given by an antisymmetric matrix $\theta_{\mu,\nu}$
and the commutation relations are given by $[x_\mu,x_\nu] = \theta_{\mu,\nu}$
can be treated as special case of this framework by multiplying
$\theta_{\mu,\nu}$ by a central element $c$. Then one calculates
the star product and obtains the classical formula, after setting back
$c$ to $1$. In the classical case also one has the formula
$f\star g = mF(f\otimes g)$ where $F\in H\otimes H$ 
is a Drinfeld twist where $H$ is the universal enveloping 
algebra of Lie algebra of formal vector fields. One would
like to have this property in general. 
Our ``normally ordered exponential'' formula for the star product 
should be rewritten in form $f\star g = mF(f\otimes g)$ where
$F$ is indeed in $H\otimes H$. In one case of interest
(``kappa-deformed space''~\cite{covKappadef,kappaind,MS})
the answer is known for all orderings. The unique choice of
an element in $H\otimes H$ is made there by trying to write our
normally ordered exponential using in addition to $\partial_mu\in H$
which may be considered as momenta operators also
some special operators which have the role of angular momenta
(defined in \cite{MS}). We hope some similar principles will enable
us to find Drinfeld twists which yield our star products in
many new cases.   

\nxpoint {\bf (Formal arguments)}
By the Hausdorff formula, using the notation from~(\ref{eq:Kfcn}),
\begin{equation*}\begin{array}{lcl}
\xi(\exp(ikx))\xi(\exp(iqx)) & = &
\exp(iK(k)\hx)\exp(iK(q)\hx) \\ &=&
\exp(iD(K(k),K(q))\hx) \\
& = &\xi(\exp(iK^{-1}(D(K(k),K(q)))x))
\end{array}\end{equation*}
where we wrote the contractions with surpressed indices.
If we denote
$$D_\phi (k,q) := K^{-1}(D(K(k),K(q))),\,\,\,\,\,\,\,\,K = K_\phi,$$
then we write this as
$
\xi(\exp(ikx))\xi(\exp(iqx)) = \xi(\exp(iD_\phi(k,q)x))
$
or equivalently
$$
\exp(ikx)\star_\phi \exp(iqx) = \exp(iD_\phi(k,q)x).
$$
In physics papers (e.g.~\cite{kappaind,MS})
$\xi(\exp(ikx))$ is usually written as
$\phi$-ordered exponential $:\exp(ik\hx):_\phi$.
Similar expressions one can write for the
deformed coproducts (in Fourier harmonics picture).
$$\begin{array}{lcl}
iD_\phi^\mu(k,q)\exp(iD_\phi(k,q)x)
& = & \partial^\mu(\exp(iD_\phi(k,q)x))\\
&=& \partial^\mu (\exp(ikx)\star_\phi \exp(iqx))\\
&=& m_\phi(\Delta_\phi(\partial^\mu)(\exp(ikx)\otimes\exp(iqx)))\\
& = & \Delta_\phi^\mu(ik,iq)(\exp(ikx)\star_\phi\exp(iqx))\\
&=& \Delta_\phi^\mu(ik,iq)\exp(iD_\phi(k,q)x)
\end{array}$$
where $\Delta_\phi^\mu(ik,iq)$ is obtained from
$\Delta_\phi(\partial_\mu)$ by substituting
$\partial^\alpha\mapsto k^\alpha$ or $q^\alpha$
depending on the tensor factor and multiplying.
Thus $iD_\phi^\mu(k,q) = \Delta_\phi^\mu(ik,iq)$.

\nxpoint Let $M_\tau := C^\lambda_{\tau\mu} x_\lambda \partial^\mu$.
The correspondence $\hat{x}_\tau\mapsto M_\tau$ is a homomorphism of
Lie algebras $\gg\to \mathrm{Lie}(A_{n,k})$ -- if we corestrict to the
image $\gg^M = \mathrm{Span}_\genfd\{M_1,\ldots,M_n\}$
and restrict the action of $\partial$-s to $\gg\subset S(\gg)$,
then this is precisely the adjoint representation. 
On the other hand, the $\gg^M\oplus\gg^*\subset A_{n,\genfd}$ is 
closed under the bracket (obviously: $[M_\tau,\partial^\rho] = -C^\rho_{\tau,\mu}\partial^\mu$, hence $\gg\cong\gg^M$
acts on $\gg^*$ here by the coadjoint representation). 

\nxpoint {\bf Theorem.} {\it Let $f\in \widehat{S(\gg^*)}$. 
Then (in symmetric ordering)

\begin{equation}\label{eq:Mlin}
M_\mu (x_\nu\star f) - x_\nu\star M_\mu f = 
M_\mu(x_\nu) f + M_\tau \chi^\tau_{\mu\nu} f
\end{equation}
where for every $1\leq \tau,\mu,\nu\leq n$, 
$\chi^\tau_{\mu\nu}\in\widehat{S(\gg^*)}$ and
$$
\chi^\tau_{\mu\nu} = \sum_{N=1}^\infty (-1)^N \frac{B_N}{N!}
\left[
C^\tau_{\mu\alpha} (\CE^{N-1})^\alpha_\nu
- (\CE^{N-1})^\tau_{\nu,\alpha}\CE^\alpha_\mu
\right].
$$
where $\CE^\alpha_\beta := C^\alpha_{\beta\rho}\partial^\rho$,
$(\CE^{N-1})^\tau_{\nu,\alpha} := 
\frac{\partial}{\partial(\partial^\alpha)}(\CE^{N-1})^\tau_\nu$,
and $M_\mu(x_\nu) = C^\lambda_{\mu\nu} x_\lambda$.
}

{\it Proof.} Write
$x_\nu\star f = x_\alpha\phi^\alpha_\nu f$ 
hence $M_\mu (x_\nu\star f) = C^\lambda_{\mu\rho}x_\lambda \partial^\rho 
x_\alpha \phi^\alpha_\nu f =  C^\lambda_{\mu\rho}x_\lambda 
(\delta^\rho_\alpha + x_\alpha \partial^\rho)\phi^\alpha_\nu f 
=  x_\nu \star M_\mu f +C^\lambda_{\mu\alpha} x_\lambda \phi^\alpha_\nu f 
- x_\alpha C^\lambda_{\mu\rho}\partial^\rho 
\phi^\alpha_{\nu,\lambda} f$, relabel the indices in the last term 
to obtain
$$\begin{array}{lcl} M_\mu (x_\nu\star f) - x_\nu\star M_\mu f &=& 
x_\tau (C^\tau_{\mu\alpha} \phi^\alpha_\nu  
- \CE^\lambda_\mu \phi^\tau_{\nu,\lambda}) f \\
& = & \sum_{N=0}^\infty (-1)^N \frac{B_N}{N!} x_\tau 
[ C^\tau_{\mu\alpha} (\CE^N)^\alpha_\nu
- (\CE^N)^\tau_{\nu,\lambda}\CE^\lambda_{\mu}]f.
\end{array}$$
For $N= 0$ only the summand 
$C^\tau_{\mu\alpha} (\CE^N)^\alpha_\nu = C^\tau_{\mu\nu}$
survives within the brackets.  For $N>1$ both summands survive, 
and within the second summand use the Leibniz rule 
for $\frac{\partial}{\partial(\partial^\lambda)}$ in the form
$(\CE^N)^\tau_{\nu,\lambda} = \CE^\tau_\rho (\CE^{N-1})^\rho_{\nu,\lambda} 
+ C^\tau_{\rho\lambda} (\CE^{N-1})^\rho_\nu$.
In the rightmost summand so obtained,
use the Jacobi identity, in the form 
$-C^\tau_{\rho\lambda} \CE^\lambda_\mu = 
-C^\tau_{\mu\lambda} \CE^\lambda_\rho+\CE^\tau_\lambda C^\lambda_{\mu\rho}$,
contracted with $(\CE^{N-1})^\rho_\nu$,
and after a cancelation of one summand, accounting for the signs 
and for the antisymmetry in lower indices, and reassembling the
$M_\tau$, one obtains the formula~(\ref{eq:Mlin}).

\end{document}